\documentclass[12pt,twoside,a4paper]{amsart}
\usepackage{amssymb}
\date{\today}

\def\A{\mathcal A}

\def\sr{\stackrel}
\def\dbar{\bar\partial}
\def\R{{\mathbb R}}
\def\C{{\mathbb C}}
\def\T{{\mathbb T}}
\def\N{{\mathbb N}}
\def\coi{{C_0^{\infty}}}
\def\cal{\mathcal}

\def\S{{\mathcal S}}
\def\B{{\mathcal B}}

\def\Im{{\rm Im\,  }}
\def\supp{{\rm supp\, }}

\def\E{{\mathcal E}}
\def\O{{\mathcal O}}
\def\L{{\mathcal L}}

\def\Re{{\rm Re\,  }}
\def\L{{\mathcal L}}

\def\ekv#1#2{\begin{equation}\label{#1}#2\end{equation}}

\def\iint{\int\hskip -2mm\int}

\def\asy{asymptotic}
\def\bdd{bounded}

\def\fu{function}

\def\fop{Fourier integral operator}

\def\hol{holomorphic}

\def\indep{independent}
\def\lhs{left hand side}

\def\neigh{neighborhood}

\def\op{operator}

\def\pb{problem}

\def\Propo{Proposition}
\def\pol{polynomial}
\def\pop{pseudodifferential operator}

\def\rhs{right hand side}
\def\sa{selfadjoint}

\def\sop{Schr{\"o}dinger operator}

\def\tf{transformation}

\def\tf{transform}

\def\ufly{uniformly}

\def\Re{{\rm Re\,}}
\def\Im{{\rm Im\,}}

\newtheorem{thm}{Theorem}[section]
\newtheorem{lma}[thm]{Lemma}

\newtheorem{prop}[thm]{Proposition}

\theoremstyle{definition}

\newtheorem{df}{Definition}

\theoremstyle{remark}

\newtheorem{preremark}{Remark}
\newtheorem{preex}{Example}

\newenvironment{remark}{\begin{preremark}}{\qed\end{preremark}}
\newenvironment{ex}{\begin{preex}}{\qed\end{preex}}

\numberwithin{equation}{section}

\begin{document}

\title[Functional calculus for non-commuting operators \dots]
{Functional calculus for non-commuting operators
with real spectra
via an iterated Cauchy formula}


\author{Mats Andersson \& Johannes Sj{\"o}strand}

\address{Department of Mathematics\\Chalmers University of Technology and
the
University of G{\"o}teborg\\SE-412 96 G{\"O}TEBORG\\SWEDEN   \&
Centre de Math{\'e}matiques\\ Ecole Polytechnique\\
FR-91128 Palaiseau cedex\\
France\\ (UMR 7640, CNRS)}

\email{matsa@math.chalmers.se \& Johannes@math.polytechnique.fr}


\thanks{First author partially supported by the Swedish Research Council,
Second author invited to G{\"o}teborg University and
  Chalmers}

 


\begin{abstract}
We define a smooth functional calculus for
a non-commuting tuple of (unbounded) operators
$A_j$ on a Banach space
with real spectra and resolvents with temperate
growth,  by means of an iterated Cauchy formula.
The construction is also extended to
tuples of more general operators allowing
smooth functional calculii.
We also discuss the relation to the
case with commuting operators.

\end{abstract}

\maketitle

\section{Introduction}\label{intro}

There are many different approaches to functional calculus for one or
several operators acting on a Banach space, a common idea being that
in order to define $f(P)$ where $P$ is some operator and $f$ a
function of some suitable class, we represent $f(x)$ as a
superposition of simpler functions $\omega _\alpha (x) $, for which
$\omega _\alpha (P)$ can be defined and then define $f(P)$ as the
corresponding superposition of the operators $\omega _\alpha (P)$. For
instance, if $P$ is a self-adjoint operator on a Hilbert space, we
have
\begin{equation}
\label{intro.1}
f(P)={1\over 2\pi }\int \widehat{f}(t)e^{itP}dt,
\end{equation}
corresponding to the representation of $f$ as a superposition of
exponential functions via Fourier's inversion formula. (Here
$\widehat{f}$ denotes the standard Fourier transform of $f$. This
approach
has been developed by M.\  Taylor \cite{Mike} and others.) Another
example is when $P$ is a bounded operator and $f$ is holomorphic in a
neighborhood of the spectrum, $\sigma (P)$, of $P$. Then
\begin{equation}
\label{intro.2}
f(P)={1\over 2\pi i}\int_\gamma f(z)(z-P)^{-1}dz
\end{equation}
where $\gamma $ is closed contour around $\sigma (P)$.

\par For problems of spectral asymptotics and scattering for partial
differential operators, the representation (\ref{intro.1}) often has
led to the sharpest known results (see H{\"o}rmander \cite{H}, Ivrii
\cite{Iv}), but the price to pay is that one has to get a good
understanding of the associated unitary group for instance via the
theory of \fop{}s or via propagation estimates. Often a formula like
(\ref{intro.2}) is easier and more practical to use. (See for
instance Agmon--Kannai \cite{Ag}, Seeley \cite{Se}.) The advantage is
that the
resolvent $(z-P)^{-1}$ can be treated with simple means (like
the theory of pseudodifferential operators).

\par If $P$ is bounded, $f(z)$ is defined with its derivatives on the
spectrum
of $P$ and has an
extension $\widetilde{f}$ to a neighborhood of the spectrum such that
$\overline{\partial }\widetilde{f}$ vanishes to infinite order on
$\sigma (P)$, and if the resolvent only blows up polynomially when $z$
tends to the spectrum, then Dynkin \cite{D} used the Cauchy-Green
formula
$$\widetilde{f}(w)=-{1\over \pi }\int (z-w)^{-1}
\partial _{\overline{z}}\widetilde{f}(z)L(dz),\ L(dz)=d(\Re z)d(\Im z),$$
to define
\begin{equation}
\label{intro.3}
f(P)=-{1\over \pi } \int (z-P)^{-1}\partial _{\overline{z}}
\widetilde{f}(z)L(dz),
\end{equation}  
and he studied the corresponding functional calculus (also with other
classes of functions $f$ allowing for wilder resolvent
behaviour). This work has been very influencial (see below).

\par Unknowingly of \cite{D}, Helffer and the second author
\cite{HeSj} used (\ref{intro.3}) as a practical device in the study of
magnetic
Schr{\"o}dinger operators in the framework of unbounded non-selfadjoint
operators $P$; $\widetilde{f}$ is then the standard almost holomorphic
extension of $f\in C_0^\infty (\R )$. (We refrain from reviewing here
the history of almost holomorphic extensions with roots in the
work of H{\"o}rmander, Nirenberg, Dynkin and others.) It was soon realized
that (\ref{intro.3}) is of great practical usefulness for many
problems in spectral and scattering theory and in mathematical
physics, because it is simple to manipulate without requiring
holomorphy of the test-functions $f$. For instance, if $P$ is an
elliptic differential operator and  $f$ belongs to a suitable class of
functions, it is very easy to show that
$f(P)$ is a pseudodifferential operator (\cite{HeSj}, \cite{DiSj}), and
other applications were obtained in cases where $f$ does not necessarily
have compact support (E.B.\ Davies \cite{Da}, A.\ Jensen, S.\ Nakamura
\cite{JeNa}). Another
application of (\ref{intro.3}) is in the area of trace formulae and
effective
Hamiltonians: For a given \op{} $P:{\cal H}\to{\cal H}$, one sometimes
introduces an auxiliary (so called Grushin-, or in more special
situations Feschbach-) problem:
\ekv{intro.4}{(P-z)u+R_-u_-=v,\ R_+v=v_+.}
Here the auxiliary operators $R_+:{\cal H}\to{\cal C_+}$, $R_-:{\cal
C}_-\to {\cal H}$ should be chosen in such a way that the problem
(\ref{intro.4}) has a unique solution
$${\cal H}\ni u=Ev+E_+v_+,\ {\cal C}_-\ni u_-=E_-v+E_{-+}v_+.$$
for all $v\in{\cal H}$, $v_+\in{\cal C}_+$. Then it is well-known that
the \op{} $E_{-+}$ inherits many of the properties of $P$, and typically
one looks for spaces ${\cal C}_\pm$ which are "smaller" in some sense, so
that the study of $E_{-+}$ may be easier than that of $P$. For trace
formulae one can show under quite general assumptions that

\ekv{intro.5}
{{\rm tr\,}f(P)={\rm tr\,}{1\over \pi }\int {\partial \widetilde{f}\over
\partial z}(z) (E_{-+})^{-1}{dE_{-+}\over dz}(z)L(dz).}
which is very useful for instance when the spaces ${\cal C}_{\pm}$
are of finite (and here equal) dimensions.

\par The approach of Dynkin \cite{D} has had a great influence on many
later works devoted to general \pb{}s of functional calculus.
In \cite{T1} J.\ Taylor
introduced a notion of joint spectrum $\sigma(P)\subset\C^n$
for
several commuting bounded operators $P_1,\ldots, P_m$ on a Banach space,
defined in terms of the mapping properties of
the operators.  This spectrum is
in general strictly smaller than the joint spectrum one obtains
by regarding   $P_j$ as elements  in some Banach subalgebra of $\L(B)$.
In \cite{T2} he then  constructed a general holomorphic functional
calculus $\O(\sigma(P))\to \L(B)$ and proved basic functorial properties.
In simple cases, for instance if the function $f$ is entire, one can
use a simple multiple Cauchy formula to represent $f(P)$, but the general
case is intricate, and Taylor's first construction was based
on quite abstract  Cauchy--Weil formulas;  later on
in \cite{T3} he made the whole construction with cohomological methods.
In \cite{A1} was given a construction based on
a multivariable notion of resolvent $\omega_{z-P}$ which permits a
representation of the
calculus analogous  to formula \eqref{intro.2}.
In special cases, for instance when the spectrum is real, such a
representation
was known earlier,
and was  used by Droste \cite{Dr}, following
Dynkin's approach \eqref{intro.3}, to obtain a smooth
functional calculus in  the multivariable case for operators with real
spectra.
This approach is extended to more general spectra in \cite{SS}.

\par Various versions of functional calculus have been used in the study
of
the joint spectrum of several commuting \sa{} operators (\cite{Co},
\cite{Ch1, Ch2}), and for nonselfadjoint operators with real spectra
in \cite{AB2}.

\par The case of non-commuting \op{}s is more difficult and more
challenging. The monograph of Nazaikinski, Shatalov and Sternin
\cite{NSS}
gives a nice treatment of such a theory and contains references to many
earlier works of V.P.\  Maslov and others. The authors build the theory on
the approximation of functions of several variables by linear
combinations
of tensor products. If $f(x_1,x_2,...,x_m)=\prod_1^m f_j(x_j)$ is such
a tensor product and $P_j$ are \op{}s on the same Banach space, that do
not
necessarily commute, it is natural to define $f(P_1,...,P_m)$ as
$f_1(P_1)\circ ...\circ f_m(P_m)$, and then approximate a general
$f(x_1,...,x_m)$ by linear combinations of tensor products, and define
$f(P_1,...,P_m)$ as the corresponding limit in the space of \op{}s. A
prototype
for non-commutative functional calculus is given by the
theory of \pop{}s, with $x_1,x_2,...,x_n,D_{x_1},...,D_{x_n}$ as the
basic
set of non-commuting operators.

\par Most approaches to
the theory of \pop{}s use direct methods rather than approximation by
tensor products. In this paper we shall suggest
a direct approach to smooth non-commutative functional calculus,
based on a multivariable version (\ref{def}) of (\ref{intro.3}).
(Another possibility, that will not be explored here is to extend
(\ref{intro.1}) 
to the multivariable case. Then, under suitable extra assumptions, one
could also consider the Weyl quantization
$$f^w(P_1,\dots ,P_m)=(2\pi
)^{-1}\int \widehat{f}(t)e^{it\cdot P}dt,$$
with $t\cdot P=\sum t_jP_j$.)
When $P_1,\dots ,P_m$ are \pop{}s with real principal symbols and $f$
belongs to 
a suitable symbol class, it will be quite obvious from our formula that
$f(P_1,\dots ,P_m)$ is also a \pop{}, by extending the arguments from
\cite{HeSj}, \cite{DiSj}. We hope that the multivariable formula
(\ref{def}) will be a useful complement to existing multivariable
functional calculii. It might provide a more direct alternative
to some parts of the theory of in \cite{NSS}. The purpose of the present
paper is merely to establish some basis for this approach and to connect
it to the one of J.\  Taylor and others (\cite{T1, T2, AB2, A1, A2}) in the
commutative case.

\par  The plan of the paper is the following:

\par In Section \ref{SpeAHol}, we introduce some special almost \hol{}
extensions of smooth \fu{}s on the real domain.

\par In Section \ref{section3} we introduce the calculus using the
formula ({\ref{def}}). and in Section \ref{further} we establish some
additional properties. Thus we get a $C_0^\infty $-calculus of several
un\bdd{} and non-commuting operators whose spectra are real and which
have
locally temperate growth of the resolvent near the real axis.

\par In Section \ref{itera}, we relate our approach to a naive iterative
approach, which amounts to treat the calculus as an \op{} valued
distribution equal to a tensor product of 1-dimensional \op{}-valued
distributions.

\par In Section \ref{Cayley}, we review the Cayley \tf{} and more general
M{\"o}bius \tf{}s of \op{}s, as a tool to reduce  many questions about
un\bdd{} \op{}s to the \bdd{} case.

\par In Section \ref{CoOp} we consider the commutative case and relate
the
theory to the Taylor approach. In particular we show that the (joint)
Taylor
spectrum and the  support of our \op{}-valued distribution agree.

\par In Section \ref{nr}, we discuss what happens when the \op{}s have
non-real spectra. In some cases there is a direct extension using
formulae
like (\ref{intro.3}) and (\ref{def}), but there are also cases where such
a functional calculus can be given differently already in the case of one
\op{} (like
for instance if we have a normal \op{} on a Hilbert space). The
conclusion
is that in all cases, one can get a multi-\op{} calculus by iterating
suitable one-dimensional formulae, in a way that is well adapted to
the spectrum of each of the individual \op{}s.

\par In Section \ref{SoFuEx}, we give some simple examples, and
show in particular that the support (unlike the joint spectrum in the
commutative case) is highly unstable under small perturbations.

\par In Section \ref{ExtFC} we extend the calculus to the case of
test-\fu{}s $f$
that do not necessarily have compact support. This is of
importance in applications to differential operators and spectral theory
(see \cite{JeNa, Da}). For simplicity, in this and the two remaining
sections,
only the case of a single \op{} is considered, with the hope that the
extension
to the multi-\op{} case should be
straight forward along the lines of the previous sections.

\par In Section \ref{recover} we show how to recover a generating \op{}
from a given homomorphism from test-functions into the \bdd{} \op{}s on
some Banach space. In the case of real spectrum it is important to have
test-\fu{}s with a non-trivial behaviour near infinity, and we give an
example of a homomorphism defined on the Schwartz-space ${\cal S}(\R )$
which is not generated by any operator.

\par In Section \ref{g(f(P))=(g(f))(P)} we establish the basic
composition
result $f(g(P))=(f\circ g)(P)$ within the framework of the extended
calculus of Section  \ref{ExtFC}

\section{Special almost holomorphic extensions}\label{SpeAHol}

\begin{lma}\label{la}
Let $f\in\coi(\R^m)$. Then there is a $\tilde f\in\coi(\C^m)$ with
support
in an
arbitrarily  small neighborhood of $\supp f$  such that
\begin{equation}\label{apa}
\partial_{\bar z_j}\tilde f=\O(|\Im z_j|^\infty), \ 1\le j\le n.
\end{equation}
\end{lma}

\begin{proof}
As a first attempt we take
\begin{equation}
\check f(z)=\frac{1}{(2\pi)^n}\int e^{iz\cdot \xi}\Big(\prod_{k=1}^m\chi
\big(\left<\xi_k\right>\Im z_k\big)\Big)
\hat f(\xi)d\xi,
\end{equation}
where $\hat f\in\S(\R^m)$ is the Fourier transform of $f$,
$$
\left<\xi_k\right>=\sqrt{1+|\xi_k|^2},\ z\cdot \xi =\sum_{k=1}^m z_k\xi
_k,
$$
and $\chi\in\coi(]-1,1[)$
is equal to $1$ in a neighborhood of $0$. Notice that the exponential
factor is bounded
on the support of the integrand so $\check f\in C^\infty(\C^m)$, and by
modifying the choice
of $\chi$ we may assume that $\check f$ has its support in an arbitrarily
small tubular
neighborhood of $\R^m$.

We have
\begin{multline*}
\partial_{\bar z_j}\check f(z)=\\ \frac{1}{(2\pi)^m}\int
 e^{iz\cdot \xi}\Big(\prod_{k=1,k\neq j}^m\chi\big(\left<\xi_k\right>\Im
z_k\big)\Big)
\left<\xi_j\right>
\frac{i}{2}\chi'\big(\left<\xi_j\right>\Im z_j\big)\hat f(\xi)d\xi.
\end{multline*}
On the support of the integrand we have $\left<\xi_j\right>\sim 1/|\Im
z_j|$
and using
the rapid decay of  $\hat f$ we get \eqref{apa}.
Clearly $\check f|_{\R^m}=f$. Notice that the map
$f\mapsto \check f$ is linear,  and at least formally
it is the tensor product of the $1$-dimensional extension
maps
\begin{equation}\label{}
\coi(\R)\ni g\mapsto \check g(z)=\frac{1}{2\pi}\int
e^{iz\xi}\chi\big(\left<\xi\right>\Im z\big)
\hat g(\xi)d\xi,
\end{equation}
 cf., Section~\ref{itera} below.
It is easy to see that (for any almost holomorphic extension $\check g$)
\begin{equation}\label{ident}
\check g(z)=\O(|\Im z|^\infty)
\end{equation}
locally uniformly when $\Re z\notin\supp g$.
In fact, if $g(x)$ has the Taylor expansion $\sum_\nu a_\nu (x-x_0)^\nu$
at some point $x_0$,  then any almost holomorphic extension must have the
expansion
$\sum_\nu a_\nu (z-x_0)^\nu$ at this point.

Let $f$ have support in $I_1\times\cdots\times I_n$, where $I_j$ are
bounded
intervals. If
$J_j\subset\subset\R$ are open intervals with $I_j\subset\subset J_j$,
let
$\psi_j\in\coi(J_j)$ be equal to $1$ near $\overline{I_j}$ and consider
\begin{equation}
\tilde f(x)=\prod_1^m\psi_j(\Re z_j)\check f(z).
\end{equation}
For $\Re z_j\in\supp\psi'_j$ we have $\check f(x)=\O(|\Im z_j|^\infty)$,
so
$\partial_{\bar z_j}\tilde f=\O(|\Im z_j|^\infty)$.

In the general case we first decompose $f$ by a partition of unity into a
finite sum of
new functions $f^\nu$, where each $f^\nu$ has support in a small box
$I_1^\nu\times\cdots\times I_m^\nu$. Then we get $\tilde f^\nu$ with
support arbitrarily close to
$I_1^\nu\times\cdots\times I_m^\nu$,  and if we sum the extensions
$\tilde f^\nu$ we get an
extension of $f$ with support in an arbitrarily small neighborhood of
$\supp f$.
\end{proof}

Notice that \eqref{apa} is stronger than the usual requirement for almost
holomorphic
extensions:
\begin{equation}
\dbar\tilde f=\O(|\Im z|^\infty).
\end{equation}
Also recall
that if $\tilde f, \check f\in\coi(\C^m)$ are almost holomorphic
extensions
of the same $f\in\coi(\R^m)$, then
\begin{equation}\label{skillnad}
\tilde f-\check f=\O(|\Im z|^\infty);
\end{equation}
this is just a special case of \eqref{ident} above.

\section{The calculus}\label{section3}

Let $P_1,\ldots ,P_m\colon\B\to \B$ be densely defined closed operators
on the
complex Banach space  $\B$. We assume that each $P_j$ has real spectrum,
\begin{equation}
\sigma(P_j)\subset\R,
\end{equation}
and that the resolvents have temperate growth locally near $\R$:
\begin{multline}\label{temp} 
{\rm For\     every\   } K\subset\subset\C \ {\rm there\  are}\ C_{K,j},
N_{K,j}\ge 0\ {\rm such\   that}\ \\
\|(z-P_j)^{-1}\|\le C_{K,j}|\Im z|^{-N_{K,j}},\   z\in K\setminus\R.
\end{multline}

\begin{df}
For $f\in\coi(\R^m)$  we put
\begin{multline}\label{def}
f(P_1,\ldots,P_m)= \\
\big(-\frac{1}{\pi}\big)^m\int\cdots\int
(\partial_{\bar z_1}\cdots\partial_{\bar z_m}\tilde
f)(z-P_1)^{-1}\cdots(z_m-P_m)^{-1}
L(dz_1)\cdots L(dz_m),
\end{multline}
where $\tilde f$ is a special almost holomorphic extension as in
Lemma~\ref{la},
and $L(dz_j)$ is the Lebesgue measure on $\C\sim \R^2$.
\end{df}

\smallskip
We first check that the right hand side of \eqref{def} is a bounded
operator on $\B$ which
depends on $f$ but not on the choice of special extension $\tilde f$.
The estimates  \eqref{apa} remain valid after differentiation so we have
for every $j$
that
$$
\partial_{\bar z_1}\cdots\partial_{\bar z_m}\tilde f=\O(|\Im
z_j|^\infty),
$$
and taking geometrical means we get
\begin{equation}
\partial_{\bar z_1}\cdots\partial_{\bar z_m}\tilde f=\O(|\Im
z_1|^\infty\cdots
|\Im z_m|^\infty).
\end{equation}
Using this in \eqref{def} we see that the integral converges in the space
of bounded
operators,  and for every $K\subset\subset\R^m$ there exist constants
$C_K, N_K \ge 0$ such that
\begin{equation}\label{uppsk}
\|I(\tilde f)\|\le C_K\sum_{|\alpha|\le N_K}\sup_K|\partial^\alpha f|
\end{equation}
for every $f\in\coi(\R^m)$ with $\supp f\subset K$, where
$I(\tilde f)$ is the right hand side of \eqref{def}.

\smallskip
Let $\check f$ be another special extension of $f\in\coi(\R^m)$.
Then
\begin{multline*}
I(\tilde f)-I(\check f)=\\
\lim_{\epsilon\searrow 0}
(-\frac{1}{\pi})^m\int\cdots\int
\big(\partial_{\bar z_1}\cdots\partial_{\bar z_m}(\tilde f-\check f)
(z_1,\ldots,z_m)\big)\times \\ \big(\prod_1^m(1-\chi(\Im
z_j/\epsilon))\big)
(z_1-P_1)^{-1}\cdots(z_m-P_m)^{-1}\prod_1^mL(dz_j),
\end{multline*}
where
$\chi\in\coi(\R)$ is equal to $1$ near the origin.
Integration  by parts gives
\begin{multline*}
I(\tilde f)-I(\check f)=\\
\lim_{\epsilon\searrow 0}
(\frac{1}{2\pi i})^m\int\cdots\int
(\tilde f-\check f)
(z_1,\ldots,z_m)\big(\prod_1^m(\chi'(\Im z_j/\epsilon))\big)\\
(z_1-P_1)^{-1}\cdots(z_m-P_m)^{-1}\prod_1^m\frac{L(dz_j)}{\epsilon}.
\end{multline*}
In view of \eqref{temp} and \eqref{skillnad},  this limit is $0$ and
hence
the definition \eqref{def} is independent of the choice of $\tilde f$.

\smallskip
It follows from \eqref{uppsk} that
\begin{equation}\label{kontupp}
\|f(P_1,\ldots,P_m))\|\le C_K\sum_{|\alpha|\le N_K}\sup_K|\partial^\alpha
f|
\end{equation}
for every $f\in\coi(\R^m)$ with $\supp f\subset K$,
which means that
$$
\coi\ni f\mapsto f(P_1,\ldots,P_m)\in\L(\B)
$$
is an operator-valued distribution
on $\R^m$.
Let $\supp(P_1,\ldots,P_m)$ denote its  support;
clearly $f(P_1,\ldots,P_m)$ is welldefined  for any
smooth $f$ defined in  some neighborhood of $\supp(P_1,\ldots,P_m)$
and vanishing  in a neighborhood of infinity.

\smallskip
Next we review Feynman notation:

\smallskip
\noindent{\bf Notation}\    If $f,P_j$ are as above and
$\pi\colon\{1,\ldots,m\}\to
\{1,\ldots,m\}$ is a permutation, we put
\begin{multline}
f(\sr{\pi(1)}{P_1},\ldots,\sr{\pi(m)}{P_m})= \\
=\big(-\frac{1}{\pi}\big)^m\int\cdots\int
(\partial_{\bar z_1}\cdots\partial_{\bar z_m}\widetilde
f)(z_1,\ldots,z_m)
(z_{\pi^{-1}(m)}-P_{\pi^{-1}(m)})^{-1}
\\(z_{\pi^{-1}(m-1)}-P_{\pi^{-1}(m-1)})^{-1}\cdots
 (z_{\pi^{-1}(1)}-P_{\pi^{-1}(1)})^{-1}
\prod_1^mL(dz_i).
\end{multline}
In simpler words, this is the same as \eqref{def} except that we
rearrange the order of
the resolvents, so that we have
$$
(z_{j_m}-P_{j_m})^{-1} (z_{j_{m-1}}-P_{j_{m-1}})^{-1}\dots
(z_{j_1}-P_{j_1})^{-1},
$$
with $\pi (j_1)=1, \phi(j_2)=2,\ldots$.

\begin{ex}[Some examples]\label{examples}
\begin{multline*}
f(\sr{3}{P_1},\sr{1}{P_2},\sr{2}{P_3})=
\big(-\frac{1}{\pi}\big)^3\int\int\int(\partial_{\bar z_1}\partial_{\bar
z_2}\partial_{\bar z_3})\tilde f
(z_1,z_2,z_3)\\
(z_1-P_1)^{-1}(z_3-P_3)^{-1}(z_2-P_2)^{-1}
L(dz_1)L(dz_2)L(dz_3).
\end{multline*}

When no indices are suspended we use the usual ordering of operators as
in
compositions, so for the operator \eqref{def}
 we have
$$
f(P_1,\ldots ,P_m)=f(\sr{m}{P_1},\ldots,\sr{1}{P_m}).
$$
This  notation can also be extended to more complicated expressions. If
 $A\in\L(\B)$,
we can define
\begin{multline*}
f(\sr{3}{P_1},\sr{1}{P_2})\sr{2}{A_{}}=\\
\big(-\frac{1}{\pi}\big)^2\int\int(\partial_{\bar z_1}\partial_{\bar
z_2})
\tilde f(z_1,z_2)
(z_1-P_1)^{-1}A(z_2-P_2)^{-1}
L(dz_1)L(dz_2).
\end{multline*}
Notice that this is {\it not} an ordinary composition of
$f(\sr{3}{P_1},\sr{1}{P_2})$
and $A$, while for instance
$$
f(\sr{2}{P_1},\sr{1}{P_2})\sr{3}{A_{}}=A\circ f(\sr{2}{P_1},\sr{1}{P_2})
$$
and
$$
f(\sr{2}{P_1},\sr{3}{P_2})\sr{1}{A_{}}= f(\sr{1}{P_1},\sr{2}{P_2})\circ
A.
$$
\end{ex}

\section{Some further properties}\label{further}

\begin{prop}\label{gurka1}
Let $f\in\coi(\R^k)$, $g\in\coi(\R^{\ell})$, $m=k+\ell$, and
$P_1,\ldots,P_m$ as above. Then
\begin{equation}\label{tensor}
f(P_1,\ldots, P_k)\circ g(P_{k+1},\ldots,P_m)=(f\otimes
g)(P_1,P_2,\ldots, P_m),
\end{equation}
where $(f\otimes g)(x_1,\ldots, x_m)=f(x_1,\ldots,x_k)g(x_{k+1},\ldots,
x_m)$.
\end{prop}

\begin{proof}
It follows  directly from the definition since we can take
$(f\otimes g)\tilde{}=\tilde f\otimes \tilde g$ as the special almost
holomorphic extension of $f\otimes g$.
\end{proof}

\begin{prop}\label{gurka2}
Let $f\in\coi(\R^m)$ and $P_1,\ldots ,P_m$ as above. If $P_{k+1}=P_k$ for
some
$k\in\{1,\ldots,k-1\}$, then
\begin{equation}
f(P_1,\ldots,P_k,P_{k+1},\ldots
,P_m)=f^{(k)}(P_1,\ldots,P_k,P_{k+2},\ldots ,P_m),
\end{equation}
where $f^{(k)}\in\coi(\R^{m-1})$ is given by
$f^{(k)}(x_1,\ldots,x_k,x_{k+2},\ldots ,x_m)=
f(x_1,\ldots,x_k,x_k,x_{k+2},\ldots,x_m)$, (i.e., by restricting $f$ to
the subspace
$x_{k+1}=x_k$).
\end{prop}

\begin{proof}
For simplicity we only consider  the case $m=2$, $k=1$, so that
$P_1=P_2=:P$. Then, using the resolvent identity,
\begin{multline*}
f(\sr{2}{P_{}},\sr{1}{P_{}})=\frac{1}{\pi^2}\int\int(\partial_{\bar z_1}
\partial_{\bar z_2}
\tilde f)(z_1,z_2)(z_1-P)^{-1}(z_2-P)^{-1}L(dz_1)L(dz_2)=\\
\frac{1}{\pi^2}
\int\int(\partial_{\bar z_1}\partial_{\bar z_2}
\tilde f)(z_1,z_2)(z_2-z_1)^{-1}(z_1-P)^{-1}L(dz_1)L(dz_2)+\\
\frac{1}{\pi^2}
\int\int(\partial_{\bar z_1}\partial_{\bar z_2}
\tilde f)(z_1,z_2)(z_1-z_2)^{-1}(z_2-P)^{-1}L(dz_1)L(dz_2)=\\
-\frac{1}{\pi}\int(\partial_{\bar z_1}\tilde
f)(z_1,z_1)(z_1-P)^{-1}L(dz_1)-
\frac{1}{\pi}\int(\partial_{\bar z_2}\tilde
f)(z_2,z_2)(z_2-P)^{-1}L(dz_2)\\
=-\frac{1}{\pi}\int\partial_{\bar z}\big(\tilde
f(z,z)\big)(z-P)^{-1}L(dz),
\end{multline*}
which gives the result since $\tilde f(z,z)$ is an almost holomorphic
extension
of $f(x,x)$.
\end{proof}

\section{Definition by iteration}\label{itera}

It is possible to construct  our functional calculus
from the single operator case by iteration.
To see this we first  extend our previous construction to
vector-valued functions.
If $f\in\coi(\R^m,\B)$ we can find a special almost holomorphic extension
and define $f(P_1,\ldots, P_m)$ in the same way as before,
just being careful to put the factor
 $\partial_{\bar z_1}\cdots\partial_{\bar z_m}\tilde f$ on the right hand
side of
all the resolvents in  formula \eqref{def}. Again this definition is
independent
of the particular choice of extension, and the estimate
\eqref{kontupp} holds.  Notice that $f(P_1,\dots
,P_m)=0$ if ${\rm supp\,}(f)\cap {\rm supp\,}(P_1,\dots ,P_m)=\emptyset$,
also when $f$
is vectorvalued (where ${\rm supp\,}(P_1,\dots ,P_m)$ is the support of
our operator-valued distribution defined initially on scalarvalued
testfunctions).
For instance, if $\phi$ is scalarvalued, $u\in \B$, and
$f(x_1,\ldots,x_m)=
\phi(x_1,\ldots,x_m)u$, then
$f(P_1,\ldots,P_m)=\phi(P_1,\ldots,P_m)u$.
Moreover, if $f(x_1,\ldots,x_k,x_{k+1},\ldots,x_m)$ is $\B$-valued,
and
$$
g(x_1,\ldots,x_k)=f(x_1,\ldots,x_k,P_{k+1},\ldots,P_m)
$$
is defined
as before, for each fixed  $(x_1,\ldots,x_k)$, then $g(x_1,\ldots, x_k)$
is
a function in $\coi(\R^k,\B)$ and
$$
f(P_1,\ldots,P_m)=g(P_1,\ldots,P_k).
$$

\begin{ex}
One can define, e.g., $f(\sr{3}{P_1},\sr{1}{P_2})\sr{2}{A}$,
cf., Example~\ref{examples}, as $g(P_1)$, where
$$
g(x_1)=A\circ f(x_1,P_2).
$$
\end{ex}

\begin{remark}\label{tensrem}
 Since we use  explicit integral formulas
the necessary verifications  for the statements above are  easily
made directly.
However one can also obtain the multi-operator calculus
in a more abstract  way.
Spaces like $\coi(\R^k)$ are nuclear, and therefore 
they behave well under topological tensor products. Since
$$\coi(\R^m,\B)=\coi(\R)\hat\otimes\cdots\hat\otimes\coi(\R)\hat\otimes\B$$
it is therefore  enough to define the functional calculus
on decomposable elements
$\phi_1(x_1)\otimes\cdots\otimes\phi_m(x_m)\otimes u$,
for $u\in\B$, which is done by the single operator calculus.
\end{remark}

As an application we  can prove

\begin{prop}
If $P_1,\ldots, P_m$  are as above, then
$$
\supp(P_1\ldots,
P_m)\subset\supp(P_1,\ldots,P_k)\times\supp(P_{k+1},\ldots,P_m).
$$
\end{prop}

\begin{proof}
Let $P=P_1,\ldots,P_k$ and $Q=P_{k+1},\ldots,P_m$,
and similarily\break $(x_1,\ldots,x_m)=(x,\xi)$.
If $\phi(x,\xi)$ has support outside
$\supp(P)\times\supp(Q)$, then
$\xi\mapsto \phi(x,\xi)$  vanishes
near $\supp(Q)$  if $x$ belongs to
(a neighborhood of) $\supp(P)$. Thus
$x\mapsto \phi(x,Q)$
vanishes in a neighborhood of $\supp(P)$  and hence
$\phi(P,Q)=0$.
\end{proof}

\begin{ex}\label{enop}
For one single operator $P$,  the support coincides with
the spectrum $\sigma(P)$, i.e., the complement of the
resolvent set.
In fact,  suppose that $f\in\coi(\R)$ has support in the  resolvent
set. Then we  may assume that $\tilde f$
has support in the resolvent set as well. However, here the resolvent
$(z-P)^{-1}$ is holomorphic, and thus
$$
-\frac{1}{\pi}\int\partial_{\bar z}\tilde f(z)(z-P)^{-1}L(dz)=
-\frac{1}{\pi}\int\partial_{\bar z}(\tilde f(z)(z-P)^{-1}) L(dz)=0
$$
by Stokes' theorem. Thus
$\supp(P)\subset\sigma(P)$.
Conversely, if $\Omega$ is  an open set in the complement
of the support, then the operator-valued function
$c(z)=(z-P)^{-1}$ has a holomorphic extension across $\R$ in
$\Omega$.
In fact, if $F\in C_0^{\infty}(\Omega)$
is an almost holomorphic extension of a function $f$ in $C_0^\infty
 (\Omega\cap\R )$, 
then it is easy to see that
$$
-\frac{1}{\pi}\int
c(\zeta)\partial_{\bar\zeta}F(\zeta)\frac{L(d\zeta)}{\zeta-z}=c(z)F(z)
$$
for each $z\in\Omega\setminus\R$. For  any given point $x^0$ in
$\Omega\cap\R$ we can choose
$F$ which  is identically one in a neighborhood, and then the integral 
provides 
the holomorphic extension at $x^0$. 
One can conclude that $\Omega$ is contained in the resolvent set of $P$. 
Thus 
$\supp(P)=\sigma(P)$.
\end{ex}

\section{The Cayley transform}\label{Cayley}

In this section we shall consider closed operators on a complex Banach
space ${\mathcal B}$ that are not necessarily densely defined. For such 
operators $P$ one defines the spectrum as usual  (namely as the 
complement 
in $\C$ of the set of $z$ for which $z-P:{\cal D}(P)\to {\cal B}$ has a 
\bdd{} inverse, where ${\cal D}(P)$ is the domain, equipped with the 
graph-norm $\Vert u\Vert +\Vert Pu\Vert $) and the spectrum 
$\sigma (P)$ becomes a closed subset of the complex plane. The point 
spectrum $\sigma _p(P)\subset \sigma (P)$ is the set of $z\in\C$ such 
that $z-P$ is not injective. In this 
section we only consider operators whose spectrum is not equal to the 
whole complex plane.

\par For any closed  operator $P$ on ${\mathcal B}$, we define its
extended spectrum $\widehat{\sigma }(P)$ as $\sigma (P)$ if $P$ is
bounded and as $\sigma (P)\cup \{\infty \}$ if $P$ is not
bounded. Then $\widehat{\sigma }(P)$ is a compact subset of the
extended plane $\widehat{\C}=\C\cup\{ \infty \}$. If $\psi $ is an
automorphism of $\widehat{\C}$, a M{\"o}bius mapping, such that $\psi
^{-1}(\infty )$ is outside the point spectrum of $P$, then $\psi (P)$
is a welldefined closed operator with extended spectrum $\psi
(\widehat{\sigma }(P))$, and it is bounded, if and only if this set is
bounded, i.e., if and only if $\psi ^{-1}(\infty )$ is outside
$\widehat{\sigma }(P)$. Moreover, $\psi (P)$ is densely defined if and
only if the range of $P-\psi ^{-1}(\infty )$ is dense (excluding the
trivial case when $\psi $ maps $\infty $ to itself, in which case $P$
and $\psi (P)$ have identical domains). More precisely, ${\mathcal
D}(\psi
(P))={\mathcal R}(P-\psi ^{-1}(\infty ))$, where ${\mathcal D}$ and
${\mathcal
R}$ indicate the domain and the range respectively. A
simple way of checking these facts is to use that if $\psi
(z)=(m_{1,1}z+m_{1,2})/(m_{2,1}z+m_{2,2})$, with $\det M\ne 0$, $M=\{
m_{j,k}\} _{1\le j,k\le 2}$, then the graph of $\psi (P)$ is equal to
$M({\rm graph\,}(P))$, where $M$ acts on ${\cal B}\times {\cal B}$ in the
natural way and ${\rm graph\,}(P)=\{ (Pu,u);\, u\in{\cal D}(P)\}$.

\par In this way, any closed operator $P$ such that $\widehat{\sigma }(P)
\subset_{\ne}\widehat{\C}$ can be transformed to a bounded
operator. If $\widehat{\sigma }(P)\subset\widehat{\R}$, one can use the
automorphism
$$C(z)={z+i\over z-i},$$
which maps $\widehat{\R}$ bijectively to the unit circle $\T$ and has the
inverse
$$z=C^{-1}(w)=i{w+1\over w-1}.$$
Thus $C$ induces a 1-1 correspondence between closed operators $A$ with
real spectra and bounded operators $B$ with $\sigma (B)\subset \T$,
such that $B-1$ is injective.

\par We also have the identity
$$ \vert w\vert ^2-1=4{\Im z\over |z-i|^2},$$
which implies that
$|\Im z|\sim d(w,\T )$, for $z$ close to $\R$ (i.e. $w$ close to
$\T$) with explicitly controled non-uniformity when
$z\to\infty$ ($w\to 1$) . Furthermore, with $A$, $B$ as above, we have
\begin{equation}
\label{kaffe}
{dw\over w-B}={A-i\over z-i}{dz\over z-A},
\end{equation}
which implies that
$(w-B)^{-1}$ has temperate growth locally near $\T_0=\T\setminus \{
1\}$ if and only if $(z-A)^{-1}$ has temperate growth locally near
$\R$.

\par If this holds, we can define a functional calculus
$$C_0^\infty (\T_0)\to {\mathcal L}({\mathcal B}),\ \phi \mapsto \phi
(B),$$
as before, by the formula
$$\phi (B)=-{1\over \pi }\int \partial _{\overline{w}}\widetilde{\phi}
  (w){L(dw)\over w-B},$$
where $\widetilde{\phi }$ is an almost holomorphic extension of $\phi $
  with compact support.

\par Clearly, $\phi \in C_0^\infty (\T_0 )$ if and only if $\phi \circ
  C\in C_0^\infty (\R )$, and as one would expect, 
\begin{equation}
\label{ct}(\phi \circ C)(A)=\phi (B).
\end{equation}
To see this, just notice that, by, (\ref{kaffe})
\begin{eqnarray*}
\phi (B)=-{1\over \pi }\int \partial _{\overline{w}}\widetilde{\phi
  }(w){L(dw)
\over w-B}={1\over 2\pi i}\int \overline{\partial }_w\widetilde{\phi
  }(w)\wedge {dw\over w-B}=\\
{1\over 2\pi i}\int \overline{\partial }_z(\widetilde{\phi }\circ
  C)(z)\wedge {A-i\over z-i}{dz\over z-A}=-{1\over \pi} \int \partial
  _{\overline{z}}(\widetilde{\phi }\circ C)(z) {A-i\over z-i}{L(dz)\over
  z-A},
\end{eqnarray*}
and the last integral is equal to $\phi \circ C(A)$ by Stokes' theorem,
since
$$(1-{A-i\over z-i}){1\over z-A}={1\over z-i}$$
is holomorphic.

\section{Commuting operators}\label{CoOp}

In this section we shall see what happens if we impose the  extra
condition
that $P_1,\ldots,P_m$ commute, but let us first recall the basic
elements of Taylor's theory for commuting operators, \cite{T1} and
\cite{T2}.
If  $A_1,\ldots,A_m$ is a tuple of commuting bounded operators
on $\B$, then there is a compact set $\sigma(A)=\sigma(A_1,\ldots,A_m)$
in $\C^m$ called the joint (Taylor) spectrum.
If $A^j$ is a sequence of commuting tuples, all of which commute
mutually, such that $A^j\to A$ in operator norm, then
$\sigma(A^j)\to\sigma(A)$ in the Hausdorff sense (this is not true
in general if they do not commute!).
For each function $f$
which is holomorphic in a neighborhood of $\sigma(A)$ one can define
$f(A)$, depending continuously on $f$,
such that it coincides with the obvious definition
if $f$ is a polynomial or entire function, and such that 
$(fg)(A)=f(A)g(A)$.
Moreover, if $f=f_1,\ldots, f_n$, and
$f(A)=f_1(A),\ldots,f_n(A)$, then
the spectral mapping property holds, i.e.,
$\sigma(f(A))=f(\sigma(A))$.

Let us now suppose that the spectrum of each $A_k$ is real. By the
spectral mapping property this holds
if and only if the joint spectrum $\sigma(A)$ is contained in $\R^m$.
Moreover,
$w\in\C^n$ is outside the spectrum if and only if
there are $C_j$ in $(A)$, the closed subalgebra of $\L(\B)$ generated
by $A_1,\ldots,A_m$,  such that
$$
\sum C_j (A_j-w_j)=1.
$$
The tuple $A$ admits a continuous  extension of the real-analytic
functional calculus to a smooth one if and only this holds
for each $A_j$, and this in turn is equivalent to the fact that
the resolvent of each $A_j$ has temperate growth in the
$\Im$-direction; it is also equivalent to that
$$
\|e^{it\cdot A}\|\lesssim  \langle t\rangle^M, \quad t\in\R^m,
$$
for some $M>0$, see, e.g., \cite{AB2}.
If  $A$ admits such a smooth functional calculus
that  extends the real-analytic functional calculus
(induced in the  natural way by the holomorphic functional
calculus),
then it is unique  and the support of the corresponding
operator-valued  distribution is precisely $\sigma(A)$.
Moreover,  there is then an operator-valued form $\omega_{z-A}$
of bidegree $(m,m-1)$ in $\C^m\setminus\sigma(A)$,
representing the resolvent of $A$, with
$$
\|\omega_{z-A}\|\le C|\Im z|^M,
$$
and the  smooth functional calculus can  be represented by
\begin{equation}\label{sparv} 
f(A)=-\int \dbar_{z}\tilde f\wedge\omega_{z-A}, 
\end{equation} 
if $\tilde f$ is a standard almost holomorphic extension of $f\in
C^{\infty}(\R^n)$,
i.e., such that $|\dbar_z\tilde f|=\O(|\Im z|^\infty)$, see \cite{AB2}.

As long as $A_k$ are bounded, our functional calculus, constructed by 
means of 
\eqref{def}, is defined for any $f\in C^{\infty}(\R^m)$, and
we  claim that
it in fact coincides with \eqref{sparv}.
To see this, let us first assume that   $f$ is
the restriction of an entire function $F$. Then we can take our special
almost holomorphic extension  to be equal to $F$ in a neighborhood of
$\R^m$,
and it  then follows from the iterated Cauchy formula that \eqref{def}
gives the holomorphic functional calculus. Since the entire functions are
dense
in $C^{\infty}(\R^m)$, the claim follows.
From the representation \eqref{sparv} it immediately follows that
the support of the functional calculus, $\supp(A)$, is equal to
$\sigma(A)$.
The same statements hold if  $\R^m$ is replaced by the real torus
$\T^m$.

\smallskip
Let us now go back to our unbounded closed operators
with real spectra.
We say that two such  operators $P_1, P_2$ commute
if the  resolvents $(z_1-P_1)^{-1}$ and $(z_2-P_2)^{-1}$ commute
for all $z_1$ and $z_2$ in the  resolvent sets. This holds if and only if
the Cayley transforms $C(P_1)$ and $C(P_2)$ commute.
If $P_1$ and $P_2$  are bounded this just means that they commute
themselves.
Now let $P_1,\ldots, P_m$ be as before,
i.e., resolvents with temperate growth,  but, in addition, commuting.
It is convenient to extend our functional calculus to the algebra
 $$
\A=\coi(\R^m)\oplus(1),
$$ 
of all
smooth functions which are constant in some neighborhood of
$\infty$.

\smallskip
Observe  that if $P_j$  are commuting, then
$$
f(P_1,\ldots, P_m)=f(\sr{\sigma(1)}{P_1},\ldots,\sr{\sigma(m)}{P_m})
$$
for any permutation $\sigma$.
From Propositions~\ref{gurka1} and \ref{gurka2}
we get

\begin{prop}\label{multi}
Suppose that $P_1,\ldots, P_m$ are as above and commuting. Then
\begin{equation}
f(P_1,\ldots,P_m)g(P_1,\ldots,P_m)=(fg)(P_1,\ldots,P_m), \quad f,g\in\A.
\end{equation}
\end{prop}

Let $C(x_1,\ldots, x_m)=(C(x_1),\ldots, C(x_m))$
be the multiple Cayley transform, and suppose that
$P_j$ are commuting and have real spectra. Then each $C(P_j)$ has
spectrum contained
in $\T$ so the joint spectrum of $C(P)$ is contained
in $\T^m$.
If all $P_j$ are bounded, then $C(z)$ is holomorphic in a neighborhood
of $\sigma(P)$ and thus $\sigma(C(P))$ is contained in
$\T^m_0=(\T\setminus\{1\})^m$
by the spectral mapping theorem. By another application of the same
theorem it follows that
\begin{equation}\label{defspek}
\sigma (P)=C^{-1}\big(\sigma(C(P))\cap \T^m_0\big).
\end{equation}
When $P_j$ are unbounded and commuting let us take \eqref{defspek}
as the definition of $\sigma(P)$.

\begin{prop}\label{spl}
If $A_j$ are as above (real spectra and temperate
resolvents) and in addition commuting, then
$$
\supp(A)=\sigma(A).
$$
\end{prop}

\begin{proof}
Let $B=C(A)$.
We are to prove that
$\sigma(B)\cap\T^m_0$ is equal to the support of
\begin{equation}
\label{(i)}
C_0^\infty (\T _0^m) \to {\mathcal L}({\mathcal B}),\quad f\mapsto f(B).
\end{equation}
By repeated use of \eqref{ct} we have that
$C(\supp(A))$ is equal to the support of (\ref{(i)}), and so the
proposition will follow.

\par To begin with, we shall extend (\ref{(i)}) to a multiplicative
mapping
\begin{equation}
\label{(ii)}
{\mathcal G}(\T ^m)\to {\mathcal L}({\mathcal B}),
\end{equation}
where ${\mathcal G}(\T ^m)$ is the class of functions in $C^\infty (\T ^m
)$ that are real-analytic in a neighborhood of $\T ^m\setminus \T
_0^m$. Let $\chi _0(t)$ be a smooth function on $\T $ which is 1 in a
neighborhood of a given compact set $K\subset \T _0$ and 0 in a
neighborhood of 1. One can find an almost holomorphic extension
$\widetilde{\chi }_0$ to a complex neighborhood of $\T $ such that
$\widetilde{\chi }_0$ is $1$ in a complex neighborhood of $K$,
and $0$ in a complex neighborhood of 1. Then
$$\widetilde{\chi }(w)=1-\prod_ {j=1}^m \widetilde{\chi }_0(w_j)$$
is identically 0 in a complex neighborhood of $K^m$ and identically 1 in
a complex neighborhood of $\T ^m\setminus \T _0^m$. After multiplication
by a cutoff function (which is $1$ in a neighborhood of $\T^m$), we may
assume that $\widetilde{\chi }$ has compact support in $\C ^m$. Now
take $f\in {\mathcal G}(\T ^m)$ and let $F$ be the holomorphic
extension at $\T ^m\setminus \T_0^m$, and $\widetilde{f}_0$ a special
almost holomorphic extension near $\T_0^m$. Then
$$\widetilde{f}=\widetilde{\chi }F+(1-\widetilde{\chi })
\widetilde{f}_0$$
is a special almost holomorphic extension of $f$
 which is even holomorphic in a complex neighborhood of $\T ^m\setminus
 \T _0^m$.

\par Since we have temperate growth of the resolvents in $\T_0^m$, we
can now  define
\begin{equation}\label{(iii)}f(B)=(-\frac{1}{\pi })^m\int ...\int
\partial_{\overline{w}_1}...\partial
 _{\overline{w}_m}\widetilde{f} \frac{L(dw_1)}{w_1-B_1}\cdots
\frac{L(dw_m)}{w_m-B_m}.\end{equation}
It is readily verified as in Section \ref{section3} that the integral
  is
independent of the choice
of $\widetilde{f}$. Also the multiplicativity follows by means of
the resolvent identity as in Proposition \ref{gurka2} so we get the
homomorphism (\ref{(ii)}).

\par Clearly (\ref{(ii)}) extends to a multiplicative mapping from
functions which are $C^\infty $ in a neighborhood of the support of
(\ref{(i)}) and real analytic in a neighborhood of $\T ^m\setminus \T
_0^m$. In particular; if $w\in \T_0^m$ is outside this support, then
(\ref{(ii)}) applies to 
$$
\phi^w_j(x)=\frac{\bar w_j-\bar x_j}{|w-x|^2},
$$
and since $\sum_j\phi^w_j(B)(w_j-B_j)=I$ it follows that
$w\notin\sigma(B)$. Thus $\sigma (B)\cap\T_0^m$ is contained in the
support of (\ref{(i)}).

\par We claim that (\ref{(ii)}) coincides with the holomorphic
  functional calculus when $f$ is real-analytic on the whole of $\T
  ^m$. In fact; if $\widetilde{f}$ is an extension with compact
  support in $\C ^m$ which is holomorphic in a complex neighborhood of
  $\T^m$, then it follows from Cauchy's formula that
$$\widetilde{f}(z)=(-{1\over \pi})^m\int ...\int \partial
 
_{\overline{w}_1}...\partial_{\overline{w}_m}\widetilde{f}(w){L(dw_1)\over
  w_1-z_1}...{L(dw_m)\over w_m-z_m}$$
there. Therefore, see e.g.,  \cite{T2}, formula (\ref{(iii)}) defines
$f(B)$
  in the holomorphic functional calculus sense, and thus it coincides
  with our definition.

\begin{lma}\label{ApproxLemma}
Suppose that $f\in C^\infty (\T ^m)$ is real-analytic in
$U\subset \T ^m$. Then there are $f_\epsilon $, $0<\epsilon \le 1$,
holomorphic
in some $\epsilon $-independent neighborhood of $\T ^m$, and a complex
neighborhood $\widetilde{U}$ of $U$, such that $f_\epsilon \to f$ in
$C^\infty (\T ^m)$ and $f_\epsilon \to f$ in ${\mathcal
  O}(\widetilde{U})$.
\end{lma}

\par To prove the lemma one defines $f_\epsilon $ by means of
convolution with a Gaussian approximation of unity, and since we can
make contour deformation in a complex neighborhood of $U$, we also get
the convergence in ${\mathcal O}(\widetilde{U})$ for a suitable
$\widetilde{U}$.

\par To see that the support of (\ref{(i)}) is contained in $\sigma
(B)$, take any $\phi \in C_0^\infty (\T_0^m)$ with support outside
$\sigma (B)$. If $\phi _\epsilon $ are as in the lemma, then $\phi
_\epsilon \to \phi $ in ${\mathcal G}(\T ^m)$, so $\phi _\epsilon (B)\to
\phi (B)$. On the other hand, since $\phi _\epsilon $ are holomorphic
in a complex neighborhood of $\sigma (B)$, and $\phi _\epsilon \to 0$
there,
$\phi _\epsilon (B)\to 0$ by the continuity of the holomorphic
functional calculus so $\phi(B)=0$. Thus Proposition \ref{spl} is
proved.\end{proof}

\begin{remark}

If $B$ is a tuple of bounded operators with $\sigma(B)\subset\T^m$ there 
is an 
 operator-valued $(m,m-1)$-form 
$\omega_{w-B}$
in $\C^m\setminus\sigma(B)$
such that  
\begin{equation} 
\label{f(B)=} 
f(B)=-\int \overline{\partial}_w\widetilde{f}\wedge \omega _{w-B}, 
\end{equation} 
if $\tilde f$ coincides with the 
holomorphic function 
$f$ in a  neighborhood 
of $\sigma(B)$ and has compact support. 
If $B$ is as in the preceding proof, it is even possible to 
choose $\omega _{w-B}$ such that 
$$\Vert \omega _{w-B}\Vert\lesssim d(w,\T^m )^{-M}$$ 
uniformly on compact sets in $\T_0^m$; this follows 
since one can define such a  form $\omega_{w-B}$
as the functional calculus \eqref{(ii)} acting on 
$s\wedge(\dbar_w s)^{m-1}$, where 
$$
s=\sum\phi_j^w(x) dw_j/2\pi i. 
$$
By Lemma~\ref{ApproxLemma}, or by a direct computation, one verifies that 
\eqref{f(B)=} can be used to define the functional calculus 
\eqref{(ii)} (if $\tilde f$ is an 
almost holomorphic 
extension which is holomorphic in a neighborhood of 
$\T^m\setminus\T^m_0$) and from this 
formula it is obvious that the support of (\ref{(i)}) is contained in 
$\sigma (B)$.
\end{remark}

\begin{prop}\label{Prop7.4}
Let $A_j$ be as above (real spectra and temperate
resolvents) and in addition commuting.
If $\phi_1,\ldots,\phi_n\in\A$, then
$\phi_j(A)$ is a commuting tuple (of bounded operators)
and $\sigma(\phi(A))=\phi(\sigma(A))$.
\end{prop}

\begin{proof}
We first prove that if $f_j\in {\mathcal G}(\T^m)$, then
$f(\sigma(B))=\sigma(f(B))$.
If $w\notin f(\sigma(B))$, then
$\phi_j(x)=(\bar w_j-\bar f_j(x))/|f(x)-w|^2$ are analytic near $\sigma
(B)$,
and according to the previous proof, $\sum_j(w_j-f_j(B))\phi_j(B)=I$, and
hence $w\notin\sigma(f(B))$. Thus $\sigma(f(B))\subset f(\sigma(B))$.

 We may assume that $f$ is real.
Assume that $f(x^0)=w$ and that $w\notin\sigma(f(B))$.
Then (since $\sigma(f(B))$ is real) we can find $C_j$, by the holomorphic
functional calculus, commuting with all $B_k$, such that
$\sum_j(w_j-f_j(B))C_j=I$. However, for each $j$ we can solve
$$
f_j(x)-w_j=\sum(x_k-x_k^0)\psi_{jk}(x)
$$
with $\psi_{jk}(x)$ in ${\mathcal G}(\T ^m)$.
It follows that
$
\sum_k (B_k-x_k^0)\sum_jC_j\psi_{jk}(B)=I,
$
and hence $x^0\notin\sigma(B)$. Thus $w\notin f(\sigma(B))$.

\smallskip
We already know that $\phi_j(A)$ are bounded and commuting.
By the definition of $\sigma(A)$, \eqref{ct}, and the first part
of the proof, we have
\begin{multline*}
\sigma(\phi(A))=\sigma(\phi\circ C^{-1}(C(A)))
=\phi\circ C^{-1}(\sigma(C(A)))=\\ =
\phi(C^{-1}(\sigma(C(A))))=\phi(\sigma(A)).
\end{multline*}
\end{proof}

We shall now see that $\phi(A)$ admits a smooth functional calculus if
$\phi=\phi_1,\ldots,\phi_n\in\A$  and $A_k$ are as in 
Proposition~\ref{Prop7.4}.
From the proposition we have that
$$
\sigma(\phi(A))=\{\xi+i\eta; \  \\
(\xi,\eta)\in\sigma(\Re\phi(A),\Im\phi(A)\}.
$$
Moreover, if $g$ is smooth in a  neighborhood of $\phi(\sigma(A))$,
then $g\circ\phi\in\A$, in the sense that it coincides with an element 
in
$\A$ in a neighborhood of $\sigma(A)$; thus $g\circ\phi(A)$ 
is defined.

\begin{prop}
Let $A_k$ be as in  Proposition~\ref{Prop7.4} and let 
$\phi=\phi_1,\ldots,\phi_n\in\A$.
If $\phi$ is real then the resolvent of each $\phi_j(A)$ has temperate
growth.

If $g$ is a smooth function in a neighborhood of $\sigma(\phi(A))$, then
\begin{equation}\label{samman}
g\circ\phi(A)=g(\phi(A))
\end{equation}
holds, if the  right hand side is defined as $\tilde
g(\Re\phi(A),\Im\phi(A))$,
where $\tilde g(\xi,\eta)=g(\xi+i\eta)$.
\end{prop}

\begin{proof}
If $g(w)$ is any polynomial in $\C^n$, then $g\circ\phi\in\A$ and
\eqref{samman} holds by Proposition~\ref{multi}. However, if $g$ is
entire, $g_N$ are polynomials,  and $g_N\to g$, then
$g_N\circ\phi\to g\circ\phi$ in $\A$ and hence
\eqref{samman} holds for all entire $g$.

If   $\phi$ is real,  it follows that
$$
\|e^{\phi(A)\cdot t}\|\le C \langle t\rangle ^M,  \quad t\in\R^m,
$$
and this  implies (is actually equivalent to) that the resolvent
of each $\phi_j(A)$ has temperate growth in the $\Im z_j$-direction.
It also implies that $\phi(A)$ admits an extension
of the holomorphic functional calculus to a smooth
functional calculus, and moreover, that
$g_N(\phi(A))\to g(\phi(A))$ if
$g_N$ are entire functions (or polynomials) and
$g_N\to g$ in $C^\infty$ in a neighborhood
of $\sigma(\phi(A))$ in $\R^n$.
It follows that \eqref{samman} holds for such $g$. The case with a
complex $\phi$ follows
by considering $\Re\phi,\Im\phi$.
\end{proof}

\section{Extension to operators with nonreal spectra}\label{nr}

In this section we shall indicate an extension of the functional calculus
to \op{}s with not necessarily real spectrum.

\par
Let ${\cal E}(\widehat{\C})$ be the space of smooth \fu{}s
on  $\widehat{\C }$, or equivalently the space of smooth \fu{}s $f(z)$,
$z\in\C$ with
$f(z)=g(1/z)$, for $\vert z\vert >1$, where $g$ is smooth on the unit
disc. If $K\subset \widehat{\C}$ is closed, let ${\cal E}(K)$ be the
space
of germs of ${\cal E}(\widehat{\C})$-\fu{}s near $K$.
We say that a closed \op{} $A$ with $\widehat{\sigma
}(A)\subset_{\ne}\widehat{\C}$ admits a smooth functional calculus
\ekv{nr.1}
{T:{\cal E}(\widehat{\sigma }(A))\to {\cal L}({\cal B}),}
if $T$ is a continuous algebra homomorphism
that extends the \hol{} \fu{}al calculus ${\cal O}(\widehat{\sigma
}(A))\to {\cal L}({\cal B})$. Such a $T$ is an ${\cal
L}({\cal B})$-valued distribution with support ${\rm supp\,}(T)$
contained in $\widehat{\sigma }(A)$, and from applying $T$ to $\phi
(z)=1/(z-w)$, $w\notin {\rm supp\,}(T)$, it follows that ${\rm
supp\,}(T)=\sigma (A)$.

\par If $A$ is \bdd{}, then $\Re z$ and $\Im z$ are in ${\cal E}(\sigma
(A))$, so $\Re A$ and $\Im A$ are \bdd{} and continuous. It also follows
that they both have real spectrum, and the continuity of $T$ implies that
their resolvents have temperate growth. We claim that
\ekv{nr.2}
{
\sigma (\Re A,\Im A)=\{ (x,y);\, x+iy\in\sigma (A)\} .
}
In fact; if we define $A^*=\Re A-i\Im A$, then $\sigma (A,A^*)$ is the
image in $\C
^2$ of $\sigma (\Re A,\Im A)$ under the biholomorphic mapping
$$(\xi ,\eta )\mapsto (z,w)=(\xi +i\eta ,\xi -i\eta ),$$
by the spectral mapping property of the \hol{} \fu{}al calculus.
Therefore,
$$\sigma (A,A^*)\subset\{ (z,w)\in\C^2;\, w=\overline{z}\},$$
and since $\sigma (A)$ is the image of $\sigma (A,A^*)$ under
$(z,w)\mapsto z$, (\ref{nr.2}) follows.
It should be emphasized that such an extension $T$ of the \hol{}
\fu{}al calculus in general is not unique.

\par We now claim that the \hol{} \fu{}al calculus $$\phi \mapsto \phi
(\Re
A,\Im A)$$ has an extension to all $\phi \in {\cal E}_{\R ^2}(\sigma (\Re
A,\Im A))$, i.e., functions $\phi$ that are smooth in some \neigh{} of
$\sigma (\Re A,\Im A)$
in $\R^2$. In fact, there is a closed ${\cal L}({\cal B})$-valued
$\overline{\partial }$-closed (2,1)-form $\omega_{(\xi,\eta)-(\Re A,\Im
A)}$ in
$\C ^2\setminus \sigma (\Re A,\Im A)$ such that $\Vert
\omega_{(\xi,\eta)-(\Re A,\Im A)} \Vert $ has temperate growth when
$\Im (\xi ,\eta )\to 0$, in view of the discussion in the previous
section. If
$\Phi (\xi ,\eta )$ is an almost \hol{} extension of $\phi $ to
$\C^2$, with compact support, then
\ekv{nr.2.5}{\phi (\Re A,\Im A)=-\int_{\C ^2} \overline{\partial }_{\xi
,\eta
} \Phi \wedge \omega_{(\xi ,\eta )-(\Re A,\Im A)}}
is an absolutely convergent integral.

\par For $f\in{\cal E}(\sigma (A))$, let $\check{f}(x,y)=f(x+iy)$. This
gives rise to an isomorphism
$${\cal E}(\sigma (A))\simeq {\cal E}_{\R ^2}(\Re A,\Im A),$$
and we claim that
\ekv{nr.3}
{f(A)=\check{f}(\Re A,\Im A)}
for all $f\in {\cal E}(\sigma (A))$, where the right hand side is defined
by \eqref{nr.2.5} and the left hand side is $T(f)$. To begin with,
(\ref{nr.3}) clearly
holds if $f$ is a real-analytic \pol{}, since the  \lhs{} is
multiplicative
by assumption and the \rhs{} has the same property as part of the \hol{}
\fu{}al calculus. The general case follows by approximation. Thus we have
found a representation of $T(f)=f(A)$ as an explicit absolutely
convergent integral
over $\C ^2$
for $f\in {\cal E}(\sigma (A))$.

\par If we have (\ref{nr.1}) but $A$ is un\bdd{}, then we just apply
first an automorphism $\psi $ of $\widehat{\C}$,  that maps $A$ to
a bounded operator $\psi(A)$ and then
 express $T(f)=f(A)=f\circ\psi^{-1}(\psi(A))$ as an absolutely convergent
integral
$$
T(f)=-\int \dbar_{\xi',\eta'} 
(F\circ\psi^{-1})\wedge\omega_{(\xi',\eta')-(\Re \psi(A),\Im\psi(A))},
$$
where $F\circ\psi^{-1}$ 
is an almost holomorphic extension of
$f\circ\psi^{-1}$,  $\C^2_{\xi',\eta'}\supset\R^2_{x',y'}$ and 
$x'+iy'=\psi(x+iy)$.

\par If we have several \op{}s $A_j$ that admit smooth \fu{}al calculii,
${\cal E}(\widehat{\sigma }(A_j))\to {\cal L}(B)$, we can define
\ekv{nr.4}
{{\cal E}(\prod \widehat{\sigma }(A_j))\to {\cal L}({\cal B})}
as an iterated integral as in Section~\ref{section3}, just taking for
$f(z_1,\ldots,z_m)\in\E(\prod\hat\sigma(A_j))$, a special almost
holomorphic
extension $\check F$ to $\C^{2m}$ of
$$\check{f}(x_1,y_1,...,x_m,y_m)=f(x+iy_1,...,x_m+iy_m)$$ such that
$$\vert \overline{\partial }_{\xi _1,\eta _1}\overline{\partial }_{\xi
_2,\eta _2}...\overline{\partial }_{\xi _m,\eta _m}\check{F}(\xi ,\eta
)\vert ={\cal O}(\vert \Im (\xi _1,\eta _1)\vert ^\infty \cdots \vert \Im
(\xi
_m,\eta _m)\vert ^\infty )$$
in a \neigh{} of $\sigma (\Re A_1,\Im A_1)\times ...\times \sigma (\Re
A_m,\Im A_m)$. In case all $A_j$ are bounded we then get  the formula
\begin{multline*}
f(A_1,\ldots,A_m)= \\ \pm \int_{\xi_1,\eta_1}\cdots\int_{\xi_m,\eta_m}
\dbar_{\xi_1,\eta_1}\dbar_{\xi_2,\eta_2}\cdots\dbar_{\xi_m,\eta_m}\check
F(\xi,\eta)
\wedge \\ \omega_{(\xi_1,\eta_1)-(\Re A_1,\Im A_1)}\wedge\ldots
\wedge \omega_{(\xi_m,\eta_m)-(\Re A_m,\Im A_m)}.
\end{multline*}
For each unbounded $A_j$  we first have to make an appropriate 
transformation with
a M{\"o}bius mapping $\psi$ as described above, but we omit the general
resulting
formula.

\begin{remark}
If $T_j$ denotes the \op{} valued distrubution 
$$T_j\phi =\phi (A_j),$$
then (\ref{nr.4}) is just the tensor product
$$T_1\otimes ...\otimes T_m$$
and it could have been defined in a more abstract way; cf.
Remark~\ref{tensrem}.
\end{remark}

\section{Some further examples}\label{SoFuEx}

The following example shows that small noncommutative  perturbations
of a pair of operators can blow up the support.

\begin{ex}
Let $\B=\C^2$ and $A$ the operator given by the matrix
$$
A=
\left(\begin{array}{rrrrrrr}

  1 & 0 \cr
                                
  0 & 0\end{array} \right),
$$
then $\sigma(A)=\{0,1\}$ and hence by the spectral mapping theorem
for commuting operators
$$
\sigma(A,A)=\{ (0,0), (1,1)\}.
$$
Now let $A_\epsilon= U^{-1}_\epsilon A U_{\epsilon}$,
where
$$
U_{\epsilon}=
\left(\begin{array}{rrrrrrr}

  \cos\epsilon  & \sin\epsilon  \cr
                                
  -\sin\epsilon  & \cos\epsilon
\end{array} \right),
$$
i.e., rotation with $\epsilon$. Then clearly
$A_\epsilon\to A$ in norm when $\epsilon\to 0$.
We claim that $\supp(A,A_\epsilon)$
is the whole product set $\{0,1\}\times\{0,1\}$.
Let us show that it contains the point
 $(0,1)$. To see this, take
smooth functions $\phi_j(x_j)$ with small supports such that
$\phi_1(x_1)$ is $1$ in a neighborhood of $0$ and
$\phi_2(x_2)$ is $1$ in a neighborhood of $1$. Then
$$
\phi_2(A_\epsilon)= U^{-1}_\epsilon \phi_2(A) U_{\epsilon}=A_\epsilon,
$$
and
$$
\phi_1(A)=
\left(\begin{array}{rrrrrrr}

  0  & 0   \cr
                                
  0  & 1
\end{array} \right)
$$
A straight forward computation shows that $f(A,A_\epsilon)=
\phi_1(A)\phi_2(A_\epsilon)$ is like
$$
\left(\begin{array}{rrrrrrr}

  0  & 0   \cr
                                
  \epsilon  & \epsilon^2
\end{array} \right)
$$
\end{ex}

Let $P_j$ and $Q_j$ be tuples as before. Using that
$$
(z-P_j)^{-1}-(z-Q_j)^{-1}=
(z-P_j)^{-1}(Q_j-P_j)(z-Q_j)^{-1},
$$
it is easy to check that
$$
\|f(Q)-f(P)\|\lesssim
\|Q-P\|=\sum \|Q_j-P_j\|.
$$
Thus if $f$ has support outside the spectrum of $P$, then
$\|f(Q)\|\lesssim\|Q-P\|,
$
so even though not zero we can at least say that $f(Q)$ is small
if $Q$ is close to $P$.

\begin{ex}
If $P$ and $Q$ are bounded (or at least if $[P,Q]$ is bounded), then
$$
[(z-P)^{-1},(w-Q)^{-1}]=
(z-P)^{-1}(w-Q)^{-1}[P,Q](z-Q)^{-1}(w-P)^{-1},
$$
and from this formula we get that
$$
\|f(P,Q)-f(Q,P)\|\lesssim
\|[P,Q]\|.
$$
It also follows that
$f(P,Q)-f(Q,P)$ is compact if $[P,Q]$ is compact.
\end{ex}

\vspace{1cm}

\def\listing#1#2#3{{\sc #1}:\ {\it #2},\ #3.}

\section{Extended functional calculus.}\label{ExtFC}
\medskip Even though everything could be
reduced by means of Cayley \tf{} to the case
of a \bdd{} \op{}, we prefer a more direct
treatment. We also restrict the attention from now on, to the case of one
single operator, and hope that the extension to the case of several
operators will turn out to be straight forward.

\subsection{The function space ${\mathcal E}$.}\label{subsection1}

\par We define ${\mathcal E}(\widehat{\R})={\mathcal E}\subset C^\infty
({\R})$
to be the space of smooth \fu{}s on ${\R}$, which posess an \asy{}
expansion,
\ekv{1.1}{f(x)\sim \sum_0^\infty a_kx^{-k},\ x\to \infty ,}
with $a_k\in {\C}$, in the sense that for every $N\in{\N}$:
\ekv{1.2}{f(x)=\sum_0^Na_kx^{-k}+x^{-N-1}r_{N+1}(x),\ \vert x\vert >1,}
where $r_{N+1}(x)$
is \bdd{} with all its derivatives.\medskip

\begin{prop}\label{Prop1.1} \it A continuous \fu{} on ${\R}$
belongs to ${\mathcal E}$ iff it has a \bdd{} extension $\widetilde{f}$
to
${\C}$ with the property that ${\partial \widetilde{f}\over
\partial \overline{z}}$is \bdd{} and satisfies
\ekv{1.3}
{{\partial \widetilde{f}\over \partial \overline{z}}(z)={\mathcal
O}_{N_0,N_1}(\langle z\rangle ^{-N_1}\vert \Im z\vert ^{N_0}),\ \forall
N_0,N_1\in{\N}.}\end{prop}

\begin{proof} Assume first that $f\in{\mathcal E}$. For
$\vert x\vert >1$, we introduce $y=-1/x$, $g(y)=f(x)$, and observe
that the existence of an \asy{} expansion (\ref{1.1}), (\ref{1.2}) is
equivalent
to the fact that $g\in C^\infty (]-1,1[)$ with $a_0=g(0)$. Let
$\widetilde{g}(y)\in C^\infty (D(0,1))$ be an almost \hol{}
extension of $g$ with
\ekv{1.4}
{{\partial \widetilde{g}\over \partial \overline{y}}(y)={\mathcal
O}_N(\vert \Im y\vert ^N),\ \forall N\in{\N}.}
Consider $\widehat{f}(x)=\widetilde{g}(-1/x)$, $x\in{\C}$, $\vert
x\vert >1$. Using that
$${\partial \over \partial
\overline{y}}=\overline{{\partial x\over \partial y}}{\partial \over
\partial \overline{x}}=\big(\overline{{\partial y\over \partial
x}}\big)^{-1}{\partial \over \partial
\overline{x}}=\overline{x}^2{\partial \over \partial \overline{x}},$$
and that $\Im y=\vert x\vert ^{-2}\Im x$, we see that
$${\partial \widehat{f}\over \partial \overline{x}}={\mathcal O}_N({\vert
\Im x\vert ^N\over \vert x\vert ^{2N}}),\ \forall N.$$
In other words, $\widetilde{f}=\widehat{f}$ satisfies (\ref{1.3}) in the
region $\vert x\vert >1$, and combining this with the standard
construction in a \bdd{} region, we get the desired extension
$\widetilde{f}$.

\par Now let $f\in C({\R})$ posess a \bdd{} continuous extension
$\widetilde{f}$ which satisfies (\ref{1.3}). Put
\ekv{1.5}
{\widetilde{g}(z)=-{1\over \pi }\int {\partial \widetilde{f}\over
\partial \overline{w}}(w) (w-z)^{-1}L(dw),}
and notice that the integral converges and that $\widetilde{g}$ is a
\bdd{} \fu{} which satisfies
$${\partial \widetilde{g}\over \partial \overline{z}}={\partial
\widetilde{f}\over \partial \overline{z}}.$$
Consequently, $\widetilde{f}-\widetilde{g}$ is a \bdd{} entire \fu{}
on ${\C}$ and hence a constant, so
\ekv{1.6}
{\widetilde{f}(z)=a_0+\widetilde{g}(z),\ a_0\in{\C}.}

\par So far we only used that
\ekv{1.7}
{{\partial \widetilde{f}\over \partial \overline{z}}(z)={\mathcal
O}(\langle z\rangle ^{-1-\epsilon }),}
for some $\epsilon >0$, and under this weaker assumption, we see that
$\widetilde{g}$ is continuous and $\widetilde{g}(z)\to 0$, $\vert
z\vert \to \infty $.

\par Now we use the full strength of (\ref{1.3}), and write
\ekv{1.8}
{{1\over w-z}=-\sum _0^{N-1} {w^k\over z^{k+1}}+{w^N\over z^N(w-z)}.}
Using this in (\ref{1.5}), we get
\begin{multline}\label{1.9}
\widetilde{g}(z)=\sum_1^N{1\over z^k}{1\over \pi }\int {\partial
\widetilde{f}\over \partial \overline{w}}(w)w^{k-1}L(dw)\\ +{1\over
z^N}(-{1\over \pi })\int {\partial \widetilde{f}(w)\over \partial
\overline{w}}w^N{1\over w-z}L(dw)
=\sum_1^N z^{-k}a_k +{1\over z^N} r_N(z)\end{multline}
with the obvious definition of $a_k$, $r_N$. Using (\ref{1.3}), we see
that
${{r_N}_\vert}_{{\R}}$is smooth and \bdd{} together with all its
derivatives. This and (\ref{1.6}) imply that $f\in{\mathcal
E}$.\end{proof}

\par Let ${\mathcal G}$ be the space of functions
$f\in {\mathcal E}$ for which the series in (\ref{1.1}) converges and is
equal
to $f(x)$ for $\vert x\vert $ sufficiently large. In other words,
${\mathcal
G}$ is the space of smooth \fu{}s on $\R$ with a
\bdd{} \hol{} extension to a domain $\{ z\in\C ; \vert z\vert >R\}$ for
some $R>0$.

\begin{prop}\label{Prop1.1.5}
A continuous \fu{} $f$ on $\R$ belongs to ${\mathcal G}$ iff it has a
\bdd{} extension $\widetilde{f}$ to $\C$, such that ${\partial
\widetilde{f}\over \partial \overline{z}}$ has compact support and
satisfies
\ekv{1.9.5}{{\partial \widetilde{f}\over \partial \overline{z}}={\mathcal
O}(\vert \Im z\vert ^N),\ \forall N\in\N .}
\end{prop}
The proof is just a slight variation of the one of \Propo{} \ref{Prop1.1}
and will be omitted.
\subsection{The \op{}}\label{subsection2}

\par Let ${\mathcal B}$ be a complex Banach space and $P:{\mathcal B}\to
{\mathcal
B}$ a densely defined closed \op{}. We assume,
\ekv{2.1}
{\sigma (P)\subset {\R},}
so that $(z-P)^{-1}\in{\mathcal L}({\mathcal B})$ is well-defined and
depends
\hol{}ally on $z\in{\C}\setminus {\R}$. Assume,
\ekv{2.2}
{\Vert (z-P)^{-1}\Vert \le {\mathcal O}(\vert \Im z\vert ^{-N_0^0}\langle
z\rangle ^{N_1^0}),}
for some fixed $N_0^0,N_1^0\in{\R}$.

\par For the ${\mathcal G}$-calculus, we will replace (\ref{2.2}) by the
weaker
assumption (\ref{temp}) (with $P=P_1$).

\subsection{The calculus}\label{subsection3}
\par For $f\in{\mathcal E}$ as in (\ref{1.1}), we recall that we have
(\ref{1.6})
where
$\widetilde{g}$ is given by (\ref{1.5}). If $P:{\mathcal B}\to
{\mathcal B}$ satisfies (\ref{2.1}), (\ref{2.2}), we define,
\ekv{3.1}{f(P)=a_01-{1\over \pi }\int {\partial \widetilde{f}\over
\partial \overline{z}}(z) (z-P)^{-1}L(dz).}
In view of (\ref{1.3}), (\ref{2.2}), this clearly defines a \bdd{}
\op{}, but we
need to check that the \rhs{} of (\ref{3.1}) only depends on $f$ and not
on
the choice of \bdd{} extension $\widetilde{f}$ satisfying (\ref{1.3}).
Let
$\check{f}$ be a second extension of $f$ with the same properties.
Then it is a standard fact that (\ref{skillnad}) holds for the difference
of the two extensions, and this estimate can also be applied to the
difference $\widetilde{g}(w)-\check{g}(w)$, where
$\widetilde{g}(w)=\widetilde{f}(-1/w)$,
$\check{g}(w)=\check{f}(-1/w)$, $\vert w\vert <1$. We conclude that for
all
$N_0,N_1\in{\N}$,
\ekv{3.3}
{(\widetilde{f}-\check{f})(z)={\mathcal O}_{N_0,N_1}(\vert \Im z\vert
^{N_0}\langle z\rangle ^{-N_1}),}
for $z\in{\C}$. From this fact and (\ref{2.2}), it is
easy to see as in Section \ref{section3}, that
$$-{1\over \pi }\int {\partial \over \partial
\overline{z}}(\widetilde{f}-\check{f})(z)(z-P)^{-1}L(dz)=0,$$
so the definition (\ref{3.1}) is indeed \indep{} of the choice of
$\widetilde{f}$.

\par Notice that the map ${\mathcal E}\ni f\mapsto f(P)\in{\mathcal
L}({\mathcal
B})$ is linear and continuous. (${\mathcal E}$ is a Frechet space with
$C^\infty $-topology for the restriction of $f\in{\mathcal E}$ to any
\bdd{} interval and the $C^\infty (]-1,1[)$-topology for the function
$f(-1/y)$.)\medskip
\begin{ex} if $\zeta \in{\ C}\setminus{\R}$,
then $(\zeta -\cdot )^{-1}\in{\mathcal E}$ and $((\zeta -\cdot
)^{-1})(P)=(\zeta -P)^{-1}$ is the resolvent.\end{ex}

\par Let us establish a basic calculus result:
\medskip
\begin{prop}\label{prop3.1} \it If $f_1,f_2\in{\mathcal E}$, then
$f_1f_2\in{\mathcal E}$, and
\ekv{3.4}
{(f_1f_2)(P)=f_1(P)f_2(P).}\end{prop}

\begin{proof} \rm write $f_j=a_{0,j}+g_j$, with $g_j(x)\sim
a_{1,j}x^{-1}+a_{2,j}x^{-2}+...$, and recall that
\ekv{3.5}
{\widetilde{g}_j(z)=-{1\over \pi }\int {\partial \widetilde{f}_j\over
\partial \overline{w}}(w)(w-z)^{-1}L(dw),\ {\partial
\widetilde{f}_j\over \partial \overline{w}}={\partial
\widetilde{g}_j\over
\partial \overline{w}}.}
Then,
\begin{multline*}f_1(P)f_2(P)=(a_{0,1}+g_1(P))(a_{0,2}+g_2(P))
\\=
(a_{0,1}a_{0,2}+a_{0,1}g_2+g_1a_{0,2})(P)+g_1(P)g_2(P),\end{multline*}
so it suffices to show (\ref{3.4}) with $f_j$ replaced by $g_j$.
This verification can be done as in the proof of Proposition \ref{gurka2}
and we omit the details.
\end{proof}
\par\noindent \it Application. \rm If $f\in{\mathcal E}$, then $\sigma
(f(P))\subset\overline{f({\R})}$, and if $\zeta \in{\C}\setminus
\overline{f({\R})}$, then
$$(\zeta -f(P))^{-1}=\big( {1\over \zeta -f}\big) (P).$$
\medskip
\par Now consider the ${\mathcal G}$-calculus and let $P$ satisfy
(\ref{2.1}), (\ref{temp}). If $f\in{\mathcal G}$, we still define $f(P)$
by
(\ref{3.1}) and show that it does not depend on the choice of
$\widetilde{f}$ as in Proposition \ref{Prop1.1.5}. Proposition
\ref{prop3.1} remains valid for the ${\mathcal G}$-calculus, and so does
the
application.

\subsection{Relation to the Cayley \tf{}.}\label{subsection4}\par
Consider the Cayley(-M{\"o}bius) transform $C$ of Section \ref{Cayley}.

\par If $f:{\R}\to {\C}$, $g:{\T}\to {\C}$,
are related by
\ekv{4.7}
{f=g\circ C,}
then $f\in{\mathcal E}={\cal E}(\widehat{\R})$ iff $g\in {\cal
E}(\T )=C^\infty ({\T})$. Let $f\in{\mathcal E}$,
$g\in C^\infty ({\T})$ be related by (\ref{4.7}).

\par With $P$ as before, define $Q\in{\mathcal L}({\mathcal B})$, by
\ekv{4.8}
{Q=C(P),}
where the \rhs{} can either be defined by our calculus or more
directly (but equivalently) as
$$C(P)=(P+i)(P-i)^{-1}=1+2i(P-i)^{-1}.$$
We know that $\sigma (C(P))\subset{\T}$, and as in Section \ref{Cayley}
we get
\ekv{4.10}
{g(Q)=f(P),}
where $G(Q)$ is defined as prior to (\ref{ct}).

\par We have the same results for the ${\cal G}$-calculus. (If $f\in{\cal
G}={\cal G}(\R )$, then $g$ belongs to the space ${\cal G}(\T )$ of
$C^\infty $-\fu{}s on $\T$ that are analytic near 1.)

\section{Recovering $P$ from the functional calculus} \label{recover}

\par In this section we show that every functional calculus ${\mathcal
E}\ni f\mapsto {\rm Op\,}(f)\in {\mathcal L}(B)$ with suitable
properties,
is of the form
${\rm Op\,}(f)=f(P)$ for some \op{} $P$ as above. We will also get the
corresponding result for the ${\mathcal G}$-calculus.

\par Assume we have a continuous linear map
\ekv{5.1}{{\mathcal E}\ni f\mapsto {\rm Op\,}(f)\in{\mathcal L}(B),}
with the property
\ekv{5.2}
{{\rm Op\,}(f_1){\rm Op\,}(f_2)={\rm Op\,}(f_1f_2),\ f_j\in{\mathcal E}.}
We further assume,
\ekv{5.3}
{\sum_{g\in C_0^\infty ({\R})}{\mathcal R}({\rm Op\,}(g))\hbox{ is
dense in }{\mathcal B},}
\ekv{5.4}
{\bigcap_{g\in C_0^\infty }{\mathcal N}({\rm Op\,}(g))=0,}
where ${\mathcal N}=$"nullspace of", ${\mathcal R}=$"range of".
\medskip
\begin{lma}\label{Lemma5.1} If $g_0\in{\mathcal E}$ satisfies
$g_0(x)\ne 0$ for all $x\in{\R}$, then ${\rm Op\,}(g_0)$ is
injective with dense range.\end{lma}

\begin{proof} If $g\in C_0^\infty $, then
$k=g/g_0\in C_0^\infty $, $g=kg_0$, so
$${\rm Op\,}(g)={\rm Op\,}(g_0){\rm Op\,}(k)={\rm Op\,}(k){\rm
Op\,}(g_0).$$
Hence
$${\mathcal R}({\rm Op\,}(g))\subset {\mathcal R}({\rm
Op\,}(g_0)),\
{\mathcal N}({\rm Op\,}(g))\supset {\mathcal N}({\rm Op\,}(g_0)),$$
and the lemma follows.\end{proof}

\par Put $\omega _z(x)=1/(z-x)$, so that $\omega _z\in{\mathcal E}$ for
$z\in{\C}\setminus{\R}$.\medskip

\begin{lma}\label{Lemma5.2}${\mathcal D}:={\mathcal R}({\rm Op\,}(\omega
_z))$, $z\in{\C}\setminus{\R}$ is \indep{} of the choice of $z$.
\end{lma}

\begin{proof} Let $z,w\in{\C}\setminus{\R}$, so
that
$\omega _w/\omega _z,\, \omega _z/\omega _w\in{\mathcal E}$. The lemma
follows from applying ${\rm Op\,}$ to the relations
$$\omega _z={\omega _z\over \omega _w}\omega _w,\ \omega _w={\omega
_w\over \omega _z}\omega _z.$$
\par \end{proof}
\medskip
\begin{df}. For $u={\rm Op\,}(\omega
_z)v\in{\mathcal D}$, $z\in{\C}\setminus{\R}$, $v\in{\mathcal B}$, we
put $Pu={\rm Op\,}(\omega _z(x)x)v={\rm Op\,}({\cdot \over z-\cdot
})v$.
\end{df}

\par We need to check that this definition does not depend on the
choice of $z,v$, in the representation of $u$, so assume that we also
have $u={\rm Op\,}(\omega _{\widetilde{z}})(\widetilde{v})$,
$\widetilde{z}\in{\C}\setminus {\R}$, $\widetilde{v}\in{\mathcal
B}$. Using that ${\rm Op\,}(\omega _z)$, ${\rm Op\,}(\omega
_{\widetilde{z}})$ are injective, we see that $\widetilde{v}={\rm
Op\,}(\omega _z/\omega _{\widetilde{z}})v$, and hence,
$$
{\rm Op\,}(x\omega _{\widetilde{z}}(x))\widetilde{v}={\rm
Op\,}(x\omega _{\widetilde{z}}){\rm Op\,}({\omega _z\over \omega
_{\widetilde{z}}})v=
{\rm Op\,}(x\omega _{\widetilde{z}}{\omega _z\over \omega
_{\widetilde{z}}})v={\rm Op\,}(x\omega _z)v.$$
Hence the definition of $P$ does not depend on the choice of $z,v$.

\par We also see that $P:{\mathcal B}\to{\mathcal B}$ is a closed \op{}
with
domain ${\mathcal D}$, with $\sigma (P)\subset{\R}$, and with
\ekv{5.5}{(z-P)^{-1}={\rm Op\,}(\omega _z).}

\par On the other hand, if $q$ is a seminorm on ${\mathcal E}$, then
\ekv{5.6}
{
q(\omega _z)\le C_0\vert \Im z\vert ^{-N_0^q}\langle z\rangle ^{N_1^q}
,}
for some $N_0^q,N_1^q\in{\N}$, and combining this with (\ref{5.5}) and
the fact that Op is continuous on ${\mathcal E}$ with values in
${\mathcal
L}({\mathcal B})$, we obtain
\ekv{5.7}
{\Vert (z-P)^{-1}\Vert \le C_0\vert \Im z\vert ^{-N_0^0}\langle
z\rangle ^{N_1^0},}
for some $N_0^0,N_1^0\in{\N}$. \medskip

\begin{prop}\label{Prop5.3} ${\rm Op\,}(f)=f(P)$ for all
$f\in{\mathcal E}$.\end{prop}

\begin{proof} \rm From (\ref{1.5}), (\ref{1.6}), we get by restriction
to the real axis,
\ekv{5.8}
{f=a_0-{1\over \pi }\int {\partial \widetilde{f}\over \partial
\overline{z}}(z)\omega _z L(dz),}
where $\widetilde{f}$ is an almost \hol{} extension of $f$ as in
\Propo{} \ref{Prop1.1}. Now (\ref{5.8}) converges in ${\mathcal E}$, so
\begin{multline*}
{\rm Op\,}(f)=a_01-{1\over \pi }\int {\partial \widetilde{f}\over
\partial \overline{z}}(z){\rm Op\,}(\omega _z)L(dz)=\\
=a_01-{1\over \pi }\int{\partial \widetilde{f}\over \partial
\overline{z}}(z)(z-P)^{-1}L(dz)=f(P),
 \end{multline*}
where we used (\ref{5.5}) for the second equality and (\ref{3.1}) for the
last
one.\end{proof}

\par ${\cal G}$ is not a Frechet space but rather an inductive limit of
such spaces:
$\lim_{R\to \infty }{\cal G}_R$, where
$${\cal G}_R=\{
f\in{\cal G};\, f\hbox{ extends to a \bdd{} \hol{} \fu{} in }\vert z\vert
>R\} .$$
A sequence of functions converges in ${\cal G}$ iff there is some $R>0$
such
that it converges in ${\cal G}_R$. Assume that we have a (sequentially)
continuous map
\ekv{5.9}
{
{\cal G}\ni f\mapsto {\rm Op}(f)\in{\cal L}({\cal B}),
}
satisfying (\ref{5.2})--(\ref{5.4}). Then we can still define a closed
densely defined operator as above. Instead of (5.7), we get (\ref{temp})
and by the same proof as above, we have
\begin{prop}\label{Prop5.4} ${\rm Op\,}(f)=f(P)$ for all
$f\in{\mathcal G}$.\end{prop}

\begin{remark}
In view of Proposition \ref{Prop5.3} it is natural to ask whether any
continuous algebra homomorphism
\begin{equation}
\label{grada}
\Phi: {\cal E}(\widehat{\R})\to {\cal L}({\cal B})
\end{equation}
corresponds to a closed operator $A$ (with a resolvent with temperate
growth as before) such that $\Phi (\phi )=\phi (A)$ for $\phi \in
{\mathcal E}(\widehat{\R})$. Given such a $\Phi $, there is a unique
homomorphism
$$\widetilde{\Phi }:\, C^\infty (\T )\to {\cal L}({\cal B}),$$
such that $\widetilde{\Phi }(f)=\Phi (f\circ C)$. If
$B=\widetilde{\Phi }({\rm id\,})=\Phi (C)$ (where ${\rm id\,}(w)=w$,
$w\in \T$), then $\widetilde{\Phi }(f)=f(B)$ for $f\in C^\infty (\T
)$, $\sigma (B)=\T $, and the resolvent has temperate growth near $\T$
(just apply to $f(z)=1/(w-z)$). If the operator $A$ exists, then
$C(A)=\Phi (C)=B$, so therefore $B-1$ must be injective. Conversely, if
$B-1$ is injective, it is easy to check that $A=C^{-1}(B)$
defines $\Phi $. (Notice that the conditions (\ref{5.3}), (\ref{5.4})
ensure that $B-1$ is injective and has dense range, respectively.)

\par The same conclusions hold if ${\cal E}(\widehat{\R})$ is replaced
by ${\cal G}(\R )$.

\par If we instead consider a similar homomorphism from ${\cal S}(\R )$
or $C_0^\infty (\R )$ things are different; then there is not
necessarily always an operator like $B$. To see this, let
$$f(x)=x(2+\sin x^m),$$
where $3\le m\in\N$ and notice that $f^*$, i.e; the composition with
$f$,
induces a continuous homomorphism ${\cal S}(\R )\to {\cal S}(\R )$. If
${\cal
  B}=H^1(\R )$, we can define a continous homomorphism ${\cal S}\to
{\cal L}(H^1({\cal \R}))$, by letting $\Phi (\phi )$ be multiplication
on $H^1(\R )$ by $f^*\phi =\phi \circ f$. It is easy to see that this
$\Phi $ cannot be extended to any function $\phi (x)=1/(z-x)$, and
therefore it does not correspond to any operator like $A$ or $B$
above.
\end{remark}

\section{A $g(f(P))=(g\circ f)(P)$ result.}\label{g(f(P))=(g(f))(P)}

\par As a preparation, we construct a suitable almost \hol{} extension
of ${\R}\ni x\mapsto (\zeta -f(x))^{-1}$, when $f\in{\mathcal E}$,
$\zeta \not\in \overline{f({\R})}$. Let $\widetilde{f}(z)$ be an
almost \hol{} extension of $f$ with
\ekv{1}
{
{\partial \widetilde{f}\over \partial \overline{z}}={\mathcal
O}_N(1)\big({\vert \Im z\vert \over \langle z\rangle ^2}\big)^N,\
\forall N\ge 0, }
and
\ekv{2}
{\nabla \widetilde{f}(z)={\mathcal O}(\langle z\rangle ^{-2}).}
Then,
\ekv{3}
{\widetilde{f}(z)=f(\Re z)+{\mathcal O}\big( {\Im z\over \langle z\rangle
^2}\big) .}
Let $\delta (\zeta )={\rm dist\,}(\zeta ,\overline{f({\R})})$. From
(\ref{3}), it follows that
\ekv{4}
{
\vert \widetilde{f}(z)-\zeta \vert >\delta (\zeta )/2,\hbox{ if
}{\vert \Im z\vert \over \langle z\rangle ^2}\ll \delta (\zeta ). }
Let $\chi \in C_0^\infty ({\R})$ be equal to 1 near 0, and put
\ekv{5}
{
\chi _\delta (z)=\chi ({C\vert \Im z\vert \over \delta \langle
z\rangle ^2}), }
where $C>0$ is large enough, but \indep{} of $\delta ,z$. Notice that
when $\delta >0$ is large enough, then $\chi _\delta (z)=1$, for all
$z\in{\C}$.

\par As an almost \hol{} extension of $x\mapsto (\zeta -f(x))^{-1}$,
we take
\ekv{6}
{
F(\zeta ,z)=\chi _{\delta (\zeta )}(z){1\over \zeta -\widetilde{f}(z)}.
}
By construction, we have
\ekv{7}
{F(\zeta ,z)={{\mathcal O}(1)\over \delta (\zeta )}.}
Further,
\ekv{8}
{{\partial \over \partial \overline{z}}F(\zeta ,z)={\partial \over
\partial \overline{z}}(\chi _{\delta (\zeta )}(z)){1\over \zeta
-\widetilde{f}(z)}+\chi _{\delta (\zeta )}(z){1\over (\zeta
-\widetilde{f}(z))^2}{\partial \widetilde{f}\over \partial
\overline{z}}(z).}
Here,
$${\partial \over \partial \overline{z}}\chi _{\delta (\zeta
)}(z)=\chi '({C\vert \Im z\vert \over \delta \langle z\rangle
^2}){\partial \over \partial \overline{z}}({C\vert \Im z\vert \over
\delta \langle z\rangle ^2})$$
has its support in a region
\ekv{9}
{
{\vert \Im z\vert \over \delta \langle z\rangle ^2}\sim 1,
}
and since
$${\partial \over \partial \overline{z}}({C\vert \Im z\vert \over
\delta \langle z\rangle ^2})={{\mathcal O}(1)\over \delta \langle
z\rangle
^2},$$
we see that the first term in the \rhs{} of (\ref{8}) is ${\mathcal
O}(\delta
^{-2}\langle z\rangle ^{-2})$ and has its support in a region (\ref{9}).
The
second term is ${\mathcal O}(1){1\over \delta ^2}({\vert \Im z\vert \over
\langle z\rangle })^N$ for all $N\ge 0$. We conclude that
\ekv{10}
{
{\partial \over \partial \overline{z}}F(\zeta ,z)={\mathcal O}_N(1)\delta
^{-2-N}({\vert \Im z\vert \over \langle z\rangle ^2})^N,\ \forall N\ge
0. }
Essentially the same estimates show that
\ekv{11}
{\nabla _zF(\zeta ,z)={\mathcal O}(1)\delta ^{-2}\langle z\rangle ^{-2}.}

\par We also notice that $${1\over \zeta -\widetilde{f}(z)}-F(\zeta
,z)=(1-\chi ({C\vert \Im z\vert \over \delta \langle z\rangle
^2})){1\over \zeta -\widetilde{f}(z)}$$
is different from 0 only when
$${\vert \Im z\vert \over \delta \langle z\rangle ^2}\ge {1\over
\widetilde{C}},$$
i.e. for
\ekv{12}{\delta (\zeta )\le {\widetilde{C}\vert \Im z\vert \over
\langle z\rangle ^2}.}

\par Now let $g$ be continuous on $\overline{f({\R})}$ with a
\bdd{} \ufly{} Lipschitz extension $\widetilde{g}(\zeta )$, $\zeta
\in{\C}$ satisfying
\ekv{13}
{
{\partial \widetilde{g}\over \partial \overline{\zeta }}={\mathcal
O}({\rm
dist\,}(\zeta ,\overline{f({\R})})^\infty ). }
Consider
$$\widetilde{h}(z)=\widetilde{g}(\widetilde{f}(z)).$$
By the chain-rule,
$${\partial \widetilde{h}\over \partial \overline{z}}={\partial
\widetilde{g}\over \partial \overline{\zeta
}}(\widetilde{f}(z))\overline{\Big( {\partial \widetilde{f}\over \partial
z}\Big) }+{\partial \widetilde{g}\over \partial \zeta }{\partial
\widetilde{f}\over \partial \overline{z}}.$$
Using that
$${\rm dist\,}(\widetilde{f}(z),f({\R}))={\mathcal O}({\vert \Im z\vert
\over \langle z\rangle ^2}),$$
and the Lipschitz properties of $\widetilde{g},\widetilde{f}$, we get
\ekv{14}
{{\partial \widetilde{h}\over \partial \overline{z}}={\mathcal
O}_N(1)\big( {\vert \Im z\vert \over \langle z\rangle ^2}\big) ^N,\
\forall N\ge 0.}
It is also clear that $\widetilde{h}$ is a \bdd{} continuous extension
of $g\circ f$ with
\ekv{15}
{
\nabla \widetilde{h}={\mathcal O}(\langle z\rangle ^{-2}).}

\par Consider
\ekv{16}
{g(f(P)):=-{1\over \pi }\int {\partial \widetilde{g}\over \partial
\overline{\zeta }}(\zeta )(\zeta -f(P))^{-1}L(d\zeta ). }
For $\zeta \in{\C}\setminus \overline{f({\R})}$, we have
\ekv{17}
{(\zeta -f(P))^{-1}=-{1\over \pi }\int {\partial \over \partial
\overline{z}}(F(\zeta ,z))(z-P)^{-1}L(dz), }
and hence,
\begin{multline}\label{18}
g(f(P))=(-{1\over \pi })^2\iint {\partial \widetilde{g}\over \partial
\overline{\zeta }}{\partial F(\zeta ,z)\over \partial
\overline{z}}(z-P)^{-1}L(dz)L(d\zeta )\\
=-{1\over \pi }\int {\partial \over \partial
\overline{z}}\Big( -{1\over
\pi }\int {\partial \widetilde{g}\over \partial \overline{\zeta }}(\zeta
)F(\zeta ,z)L(d\zeta )\Big) (z-P)^{-1}L(dz), \end{multline}
where the first double integral converges in operator norm, so the same
holds for the $\int (...) L(dz)$ integral in the last expression,
which we can view as
\ekv{19}
{
\lim_{\epsilon \to 0}-{1\over \pi }\int (1-\chi _\epsilon (z)){\partial
\over \partial
\overline{z}}(...)(z-P)^{-1}L(dz). }

\par Consider
\begin{multline*}
-{1\over \pi }\int {\partial \widetilde{g}\over \partial
\overline{\zeta }}(\zeta )F(\zeta ,z)L(d\zeta )\\
=-{1\over \pi }\int {\partial \widetilde{g}\over \partial
\overline{\zeta }}(\zeta ){1\over \zeta -\widetilde{f}(z)}L(d\zeta
)+{1\over \pi }\int {\partial \widetilde{g}\over \partial
\overline{\zeta }}(\zeta ) {1-\chi _{\delta (\zeta )}(z)\over \zeta
-\widetilde{f}(z)}L(d\zeta )\\
=\widetilde{g}(\widetilde{f}(z))+ {1\over \pi }\int {\partial
\widetilde{g}\over \partial
\overline{\zeta }}(\zeta ) {1-\chi _{\delta (\zeta )}(z)\over \zeta
-\widetilde{f}(z)}L(d\zeta ).\end{multline*}
As already observed, the integrand in the last integral is $\ne 0$
only for $\delta (\zeta )\le {\mathcal O}(1){\vert \Im z\vert \over
\langle z\rangle ^2}$, and using that ${\partial \widetilde{g}\over
\partial \overline{\zeta }}(\zeta )={\mathcal O}(\delta (\zeta )^\infty
)$, we see that
\ekv{20}
{-{1\over \pi }\int {\partial \widetilde{g}\over \partial
\overline{\zeta }}(\zeta )F(\zeta ,z)L(d\zeta
)=\widetilde{g}(\widetilde{f}(z))+{\mathcal O}(1)\big( {\vert \Im z\vert
\over
\langle z\rangle ^2}\big)^\infty .}
Using this in the last integral in (\ref{18}), represented as a limit as
in
(\ref{19}), together with the temperate growth of the resolvent, we get
\ekv{21}
{g(f(P))=-{1\over \pi }\int {\partial \over \partial
\overline{z}}(\widetilde{g}(\widetilde{f}(z)))(z-P)^{-1}L(dz)=(g\circ
f)(P).}

\end{document}